\documentclass[11 pt]{article}
\usepackage{amsfonts,amssymb}
\usepackage{lipsum,amsmath,amsthm}
\usepackage{ dsfont }
\usepackage{appendix}
\usepackage{color}
\usepackage{esint}
\usepackage{mathtools}
\usepackage{scalerel}
\usepackage{comment}
\usepackage[right = 2.5cm, left=2.5cm, top = 2.5cm, bottom =2.5cm]{geometry}
\usepackage[utf8]{inputenc}
\usepackage[english]{babel}

\usepackage{graphicx} 

\usepackage[utf8]{inputenc}
\usepackage[english]{babel}

\usepackage{color}
\usepackage{enumerate}
\usepackage{enumitem}
\usepackage{fancyhdr}
\usepackage{upgreek}
\usepackage[T1]{fontenc}

\usepackage{accents}
\usepackage{centernot}
\usepackage{stmaryrd}
\usepackage[super]{nth}
\usepackage[normalem]{ulem}

\usepackage{amsmath}
\usepackage{amscd}
\usepackage{amsfonts}
\usepackage{amsthm}
\usepackage{amstext}
\usepackage{amssymb}
\usepackage{amsbsy}
\usepackage{mathrsfs}
\usepackage{mathtools}
\usepackage{array}
\usepackage{tikz-cd}
\usepackage{bm}
 \usepackage{url}

\usepackage{scalerel}[2014/03/10]
\usepackage[usestackEOL]{stackengine}
\def\intavg{\,\ThisStyle{\ensurestackMath{%
    \stackinset{c}{0\LMpt}{c}{0\LMpt}{\SavedStyle-}{\SavedStyle\phantom{\int}}}%
    \setbox0=\hbox{$\SavedStyle\int\,$}\kern-\wd0}\int}

\def\XXint#1#2#3{{\setbox0=\hbox{$#1{#2#3}{\int}$ }
\vcenter{\hbox{$#2#3$ }}\kern-.6\wd0}}

\usepackage{float}
\usepackage{caption}
\usepackage{subcaption}
\usepackage{graphicx}

\usepackage{tikz}
\usetikzlibrary{decorations.markings}
\usetikzlibrary{decorations.pathreplacing}
\usetikzlibrary{calligraphy}
\usetikzlibrary{positioning}
\usetikzlibrary{arrows}
\usetikzlibrary{fadings}
\usetikzlibrary{math}

\usepackage[nobysame]{amsrefs}

\usepackage[linktocpage]{hyperref}
\hypersetup{
	colorlinks=false,
	linkcolor=cyan,
	citecolor=magenta,
	filecolor=cyan,    
	urlcolor=magenta,
}

\makeatletter
\renewcommand*{\@fnsymbol}[1]{\ensuremath{\ifcase#1\or \or \or \or
   *\or \dagger\or \ddagger\or \mathsection \else\@ctrerr\fi}}
\makeatother

\newtheorem{theorem}[subsection]{Theorem}

\newtheorem{lemma}[subsection]{Lemma}

\newtheorem{proposition}[subsection]{Proposition}
\newtheorem{corollary}[subsection]{Corollary}

\newtheorem*{conjecture*}{Conjecture}

\theoremstyle{definition}
\newtheorem{definition}[subsection]{Definition}

\theoremstyle{remark}
\newtheorem{remark}[subsection]{Remark}

\begin{document}
\title{A Unified Approach to Mixing and Regularity for Passive Scalar Transport by Sobolev Vector Fields\footnote{\textit{Date}: \today.}\footnote{\textit{MSC}: 35F05, 35Q35, 76B03, 76F25, 76R10}\footnote{\textit{Keywords}: Passive scalar transport; mixing estimates; Sobolev vector fields; harmonic analysis.}}
\author{Lucas Huysmans\thanks{Max Planck Institute for Mathematics in the Sciences, Inselstra{\ss}e 22, 04103 Leipzig, Germany. \hfill \break \textit{Email address}: \texttt{lucas.huysmans@mis.mpg.de}} \and Ayman Rimah Said\thanks{Chargé de recherche au  CNRS at Laboratoire de Mathématiques de Reims (LMR - UMR 9008). \hfill \break \textit{Email address}: \texttt{ayman.said@cnrs.fr}}}
\date{ \ }
\maketitle

\begin{abstract}

We develop a new framework for quantitative estimates of passive
scalar transport along Sobolev vector fields in $W^{1,p}$, when
$p>1$. Our approach is based on Christ-Journ\'{e} singular integral estimates. We show (i) a new stability estimate which quantifies
the dependence of the solution on specific frequencies of the initial data;
(ii) a new exponential mixing bound in the full DiPerna-Lions well-posedness
class; (iii) propagation of ``logarithmic'' Fourier regularity of the
solution; (iv) quantitative convergence rates for the vanishing diffusivity
and mollification limits; and (v) a logarithmic decay rate for the standard
DiPerna-Lions commutator.

\end{abstract}

\section{Introduction}
    In this paper, we present a general methodology to compute quantitative estimates for the transport of a passive scalar by a divergence-free vector field in the Sobolev space $W^{1,p}$ when $p>1$. We show that the Christ-Journ\'e singular integral estimates in \cite{seeger2019multilinear} give commutator estimates on the transport of specific frequencies of the passive scalar (Lemma \ref{thm:calderon}). By making appropriate choices of the singular kernel, we show that this harmonic approach unifies and improves many of the estimates in the current literature.

    Our main contribution is a new stability estimate for the Cauchy problem,
which quantifies the continuous dependence of the solution on specific
frequencies of the initial datum (Section~\ref{sec:weakstability}). From
this, we deduce a new exponential mixing bound for the entire DiPerna-Lions
well-posedness class (Section~\ref{sec:mixing}), as well as propagation
of mild ``logarithmic'' Sobolev regularity of the passive scalar (Section~\ref{sec:propagation}). We also prove new logarithmic convergence rates
for the vanishing diffusion and mollification limits for the passive scalar
transport equation (Section~\ref{sec:convergence}). Finally, we discuss
how the harmonic estimates in this paper are related to logarithmic decay
rates for the standard DiPerna-Lions commutator (Section~\ref{sec:diperna}).

\medskip
    
    We point the reader to the lecture notes \cite{ambrosio2014continuity} for a complete introduction to the weak formulation of the transport equation. Let $\rho(x,t):\mathbb{R}^d\times[0,T]\to\mathbb{R}$ be a weak solution to the transport equation
    \begin{equation}\label{eq:TE}
        \begin{cases}
            \frac{\partial \rho}{\partial t}(x,t)+\nabla \cdot(u(x,t)\rho(x,t))=0, \\
            \rho(x,0)=\rho_0(x),
        \end{cases}
    \end{equation}
    along a divergence-free vector field $u(x,t):\mathbb{R}^d\times[0,T]\to\mathbb{R}^d$, $\nabla\cdot u=0$, with $u(x,t)\in L_t^1 W_x^{1,p}$. Suppose $\rho(x,t)\in L_t^\infty L_x^{q}$, and denote by $\rho_0(x)\in L_x^{q}$ the initial datum, and also by $\rho_T(x)\in L_x^{q}$ the weak trace of $\rho(x,t)$ at time $t=T$. We will denote by $\frac{1}{r}=\frac{1}{p}+\frac{1}{q}$ and assume throughout that $1<p,q\le\infty$, $1\le r<\infty$.

    We remark that the assumption $p>1$ is also necessary in most previous
work, see for instance
\cite{crippa2008estimates,seis2013maximal,leger2018new,hadvzic2018singular},
since the reliance on various harmonic estimates fails when $p=1$. It is
a significant open problem to extend any of these results to the case
$p=1$ (for which well-posedness also holds, see
\cite{diperna1989ordinary,ambrosio2004transport}). This is further related
to a conjecture on the maximal rate of mixing in $W^{1,1}$ by Bressan
\cite{bressan2003lemma}. We now state our results, whose proofs are located
in the remainder of the paper.

    \medskip
    
    The key approach of this paper, inspired by work on mixing in \cite{leger2018new,hadvzic2018singular}, is that if $p,q>1$, the Christ-Journ\'e singular integral estimate in \cite{seeger2019multilinear} controls the transport of the specific frequencies $(\rho * K)(x,t)$ of the solution for any convolution kernel $K(x)$. Denoting the error in the transport of $(\rho*K)(x,t)$ by $R(x,t)$:
    \begin{equation*}
        \left(\frac{\partial}{\partial t} + u(x,t)\cdot\nabla\right)(\rho * K)(x,t) = R(x,t),
    \end{equation*}
    then we have the following bound on the commutator $R(x,t)$.
    {\renewcommand{\thesubsection}{\ref{thm:calderon}}
    \begin{lemma}
        For the commutator $R(x,t):\mathbb{R}^d\times[0,T]\to\mathbb{R}$ given by
        \begin{equation*}
            R(x,t) = u\cdot\nabla(\rho*K) - \nabla\cdot\left((u\rho)*K\right),
        \end{equation*}
        we have the bound
        \begin{equation}\label{eq:frequencycascadeintro}
            \left\|R(\cdot,t)\right\|_{L_x^r}\le C_{p,q}\left\|\rho(\cdot,t)\right\|_{L_x^q}\left\|\nabla u(\cdot,t)\right\|_{L_x^p}\left(\sum_{i,j=1}^d\left\|\frac{\partial(x_j K)}{\partial x_i}\right\|_{CZ}\right) \quad \hbox{for} \quad \frac{1}{r}=\frac{1}{p}+\frac{1}{q},
        \end{equation}
        with $1<p,q\le\infty$, $1\le r<\infty$ and $C_{p,q}>0$ depends also on the dimension $d\ge1$. $\left\|\cdot\right\|_{CZ}$ refers to the Calderon-Zygmund norm, i.e. the minimum constant
        \begin{equation*}
            \|K\|_{CZ} = \mathrm{min}\left\{\left\||x|^dK\right\|_{L_x^\infty},\left\||x|^{d+1}\nabla K\right\|_{L_x^\infty},\left\|\hat{K}\right\|_{L_\xi^\infty}\right\},
        \end{equation*}
        where $\hat{K}(\xi):\mathbb{R}^d\to\mathbb{C}$ denotes the Fourier transform of $K(x):\mathbb{R}^d \to \mathbb{R}$.
    \end{lemma}
    }

    In the remainder of this introduction, we summarise the results, the proofs
of which are located in the main body of the paper. First, we give a new
stability estimate by making an appropriate choice of the kernel
$K(x)$ (Section~\ref{sec:weakstability}). We then give an overview of how
our method allows for improving several well-studied estimates in the current
literature on the passive scalar transport equation (see Sections~\ref{sec:mixing}, \ref{sec:propagation}, \ref{sec:convergence}). Across
these sections, it demonstrates how this harmonic approach provides a more
straightforward, unified framework for doing quantitative analysis on a
series of well-studied problems.

    \setcounter{subsection}{0}
    \subsection{Quantitative well-posedness and stability of the Cauchy problem}\label{sec:weakstability}

    Our main result is a quantitative estimate on the dependence of the solution
$\rho (x,t)$ on specific frequencies of the initial datum
$\rho _{0}(x)$. By constructing a kernel in \eqref{eq:frequencycascadeintro} to control the transfer of mass between
specific length scales/frequencies, we give the first such quantification
of well-posedness of the Cauchy problem.

        We say that a mollifier $\varphi(x)\in \mathcal{S}(\mathbb{R}^d;\mathbb{R})$ (a Schwartz function with unit mass, $\int\varphi(x)dx=1$) satisfies the frequency cutoff assumption if its Fourier transform $\hat{\varphi}(\xi)\in\mathcal{S}(\mathbb{R}^d;\mathbb{C})$ satisfies the following:
        \begin{equation}\label{eq:frequencycutoffintro}
            \begin{gathered}
                \hat{\varphi}(\xi)=1 \text{ for } |\xi|\le 1, \\
                \hat{\varphi}(\xi)=0 \text{ for } |\xi|\ge 2.
            \end{gathered}
        \end{equation}
        For any $\delta>0$, we shall denote the rescaled mollifier by
        \begin{equation*}
            \varphi_\delta(x) = \delta^{-d}\varphi\left(\frac{x}{\delta}\right).
        \end{equation*}
        Our main result is then the following.

        {\renewcommand{\thesubsection}{\ref{thm:quantitative2}}
        \begin{theorem}
            Let $\rho(x,t):\mathbb{R}^d\times[0,T]\to\mathbb{R}$, with $\rho(x,t)\in L_t^\infty L_x^q$, and initial datum $\rho_0(x)\in L_x^q$, be a weak solution to the forced transport equation
            \begin{equation*}
                \begin{cases}
                    \frac{\partial \rho}{\partial t}(x,t)+\nabla \cdot(u(x,t)\rho(x,t))=f(x,t), \\
                    \rho(x,0)=\rho_0(x),
                \end{cases}
            \end{equation*}
            along a divergence-free vector field $u(x,t):\mathbb{R}^d\times[0,T]\to\mathbb{R}^d$, with $\nabla u(x,t)\in L_t^1L_x^p$, and distributional force $f(x,t)\in\mathcal{S}'(\mathbb{R}^d\times\mathbb{R};\mathbb{R})$ defined as a Schwartz distribution supported on $\mathbb{R}^d\times[0,T]$.
    
            Let $\varphi(x)\in\mathcal{S}(\mathbb{R}^d;\mathbb{R})$ be a standard mollifier satisfying the frequency-cutoff assumption \eqref{eq:frequencycutoffintro}.
    
            Then, for any $\delta>0$ and $\kappa > 0$, we have the following quantitative bound on the transfer of mass between frequencies:

            For $\frac{1}{r}=\frac{1}{p}+\frac{1}{q}$ and $1<p,q\le\infty$, $1\le r < \infty$, if
            \begin{equation}\label{eq:exponentialcascadeintro}
                \delta_0 \le \frac{1}{4}\delta \exp\left(-\frac{C_{p,q}}{\kappa}\left\|\rho\right\|_{L_t^\infty L_x^q}\left\|\nabla u\right\|_{L_t^1 L_x^p}\right),
            \end{equation}
            then
            \begin{equation}\label{eq:weakstabilityintro}
                \left\|\rho*\varphi_{\delta}\right\|_{L_t^\infty L_x^r} \le \kappa + A\left\|\rho_0*\varphi_{\delta_0}\right\|_{L_x^r} + A\left\|f*\varphi_{\delta_0}\right\|_{L_t^1 L_x^r},
            \end{equation}
            where the constant $C_{p,q} > 0$ depends only on the parameters $1<p,q\le\infty$, the choice of mollifier $\varphi(x)\in\mathcal{S}(\mathbb{R}^d;\mathbb{R})$, and the dimension $d\ge1$, while $A>0$ depends only on the choice of mollifier $\varphi(x)\in\mathcal{S}(\mathbb{R}^d;\mathbb{R})$ and the dimension $d\ge1$.
        \end{theorem}
        }



        Theorem~\ref{thm:quantitative2} is a bound on the amount of mass
$\kappa >0$ that can transfer between any two length scales
$0<\delta _{0}<\frac{1}{4}\delta $, with the scale gap \eqref{eq:exponentialcascadeintro} growing at an exponential rate in the
norm $\left \|\nabla u\right \|_{L_{t}^{1}L_{x}^{p}}$. As an immediate
corollary, we recover the standard well-posedness/continuous dependence
of $\rho (x,t)$ on $\left \|\rho _{0}\right \|_{L_{x}^{r}}$ of DiPerna-Lions
\cite{diperna1989ordinary} when $p>1$, by taking the appropriate limit
$\kappa ,\delta _{0}\to 0$ in \eqref{eq:weakstabilityintro}.


    \setcounter{subsection}{1}
    \subsection{Mixing}\label{sec:mixing}
    Fixing some suitably small $\kappa >0$, Theorem~\ref{thm:quantitative2} then says that essentially no mass can transfer
from below $\delta _{0}$ to above $\delta $, or vice versa (by time-reversibility).
Such estimates are often referred to as mixing bounds, see for instance
\cite{crippa2008estimates,alberti2019exponential}, since they control how
quickly the passive scalar can approach a homogeneous, or ``mixed'' final
state with only small scale fluctuations present. Indeed in
\cite{alberti2019exponential} the authors refer to bounds on quantities
such as
$\left \|\rho (\cdot ,t)*\varphi _{\delta}\right \|_{L_{x}^{\infty}}$ as
control on the ``geometric mixing scale''. One can therefore view Theorem~\ref{thm:quantitative2} precisely as a quantitative geometric mixing estimate.
        
        We define the negative Sobolev norm $W_x^{-1,r}$ as the dual of $\left\|\phi\right\|_{W_x^{1,r'}}:=\left\|\phi\right\|_{L_x^{r'}} + \left\|\nabla \phi\right\|_{L_x^{r'}}$, where $\frac{1}{r'}+\frac{1}{r}=1$.
        
        {\renewcommand{\thesubsection}{\ref{thm:mixing}}
        \begin{theorem}
            Let $\varphi(x)\in L^1(\mathbb{R}^d;\mathbb{R})$ be any positive mollifier, $\varphi(x)\ge0$. For any $0<\delta\le 1$, such that
            \begin{equation*}
                \left\|\rho_0*\varphi_\delta\right\|_{L_x^r}\ge\frac{1}{2}\left\|\rho_0\right\|_{L_x^r},
            \end{equation*}
            i.e. $\rho_0(x)$ is \textit{not} already mixed at scale $\delta>0$, then
            \begin{gather*}
                \left\|\rho_T\right\|_{W_x^{-1,r}}\ge \delta A\left\|\rho_0\right\|_{L_x^r}\exp\left(-C_{p,q}\frac{\left\|\rho_0\right\|_{L_x^q}}{\left\|\rho_0\right\|_{L_x^r}}\left\|\nabla u\right\|_{L_t^1L_x^p}\right),
            \end{gather*}
            for $\frac{1}{r}=\frac{1}{p}+\frac{1}{q}$ with $1<p,q\le\infty$, $1\le r<\infty$. The constants $A,C_{p,q}>0$ also depend on the dimension $d\ge1$ and the choice of positive mollifier $\varphi(x)\in L^1(\mathbb{R}^d;\mathbb{R})$.
        \end{theorem}
        }

         In Theorem~\ref{thm:mixing} we give new explicit expressions for the
mixing rate
$C_{p,q}
\frac{\left \|\rho _{0}\right \|_{L_{x}^{q}}}{\left \|\rho _{0}\right \|_{L_{x}^{r}}}$,
and initial constant
$\delta A \left \|\rho _{0}\right \|_{L_{x}^{r}}$. In particular, this
is the first mixing bound valid for the entire DiPerna-Lions well-posedness
class $\rho _{0}(x)\in L_{x}^{q}$.

Mixing estimates of this form have a long history of study; see for instance
\cite{seis2013maximal,iyer2014lower,leger2018new} for bounds of the form
%
\begin{equation}
\label{eq:intromixingbound}
\left \|\rho _{T}\right \|_{H_{x}^{-1}} \ge A(\rho _{0})\exp \left (-C(
\rho _{0})\int _{0}^{T} \left \|\nabla u(\cdot ,t)\right \|_{L_{x}^{p}}
\;dt\right ).
\end{equation}
We remark that the dependence of the mixing rate in Theorem~\ref{thm:mixing},
$C(\rho _{0})=C_{p,q}
\frac{\left \|\rho _{0}\right \|_{L_{x}^{q}}}{\left \|\rho _{0}\right \|_{L_{x}^{r}}}$,
on the initial datum is necessary by dimensional analysis, since
$\left \|\nabla u\right \|_{L_{t}^{1}L_{x}^{p}}$ has units of length to
the power $\frac{d}{p}$. The conserved ratio
$
\frac{\left \|\rho _{0}\right \|_{L_{x}^{q}}}{\left \|\rho _{0}\right \|_{L_{x}^{r}}}$
has units of length to the power $-\frac{d}{p}$, and so restores the correct
scaling. Note that when $p=\infty $ this fraction cancels, and we recover
the classical exponential rate
$\left \|\nabla u\right \|_{L_{t}^{1}L_{x}^{\infty}}$ as in the more standard
Cauchy-Lipschitz theory, see for instance
\cite{ambrosio2014continuity}. Meanwhile, the dependence of the initial
constant in Theorem~\ref{thm:mixing},
$A(\rho _{0})=\delta A\left \|\rho _{0}\right \|_{L_{x}^{r}}$, on
$\delta >0$ is precisely the ``geometric mixing scale''
\cite{alberti2019exponential} of the initial datum.

Both the expression for the mixing rate, and for the initial constant,
improve upon previous bounds
\cite{seis2013maximal,iyer2014lower,leger2018new,hadvzic2018singular},
which are not valid for the DiPerna-Lions class
$\rho _{0}(x)\in L_{x}^{q}$. More precisely, in
\cite{seis2013maximal,leger2018new,hadvzic2018singular} the authors require
additional smoothness of the initial datum, while in
\cite{iyer2014lower} one requires $\rho _{0}(x)\in L_{x}^{\infty}$ and
the constants are less explicit than Theorem~\ref{thm:mixing}.

    \setcounter{subsection}{2}
    \subsection{Propagation of regularity}\label{sec:propagation}
        We now give a different application of the commutator estimate \eqref{eq:frequencycascadeintro}. It is well-known \cite{leger2018new,ben2019convergence,brue2021sharp,brue2021advection,meyer2022propagation} that a passive scalar transported by a Sobolev vector field propagates some mild "logarithmic" regularity, for example $\|\log(1-\Delta)^\frac{p}{2}\rho\|_{L_t^\infty L_x^2}$ when $p>1$, provided the initial datum is also similarly regular.
        
        We prove that such logarithmic Fourier regularity is also propagated in the Lebesgue space $L_x^r$, by taking the kernel $K(x)=\log(1-\Delta)$ in \eqref{eq:frequencycascadeintro}:
        {\renewcommand{\thesubsection}{\ref{thm:propagation}}
        \begin{proposition}
            For the convolution kernel $K(x)=\log(1-\Delta)$, given by the Fourier multiplier
            \begin{equation*}
                \hat{K}(\xi)=\log(1+\left|\xi\right|^2),
            \end{equation*}
            one has propagation of
            \begin{equation*}
                \frac{d}{dt}\left\|\log(1-\Delta)\rho(\cdot,t)\right\|_{L_x^r}\le C_{p,q} \left\|\rho(\cdot,t)\right\|_{L_x^q}\left\|\nabla u(\cdot,t)\right\|_{L_x^p} \quad \hbox{for} \quad \frac{1}{r}=\frac{1}{p}+\frac{1}{q},
            \end{equation*}
            where the constant $C_{p,q}>0$ depends only on the parameters $1<p,q\le\infty$ and dimension $d\ge1$.
        \end{proposition}
        }

        This result aligns with the optimal Besov-type estimate
in \cite[Theorem 1.1]{brue2021advection} ($q=\infty $). However, it is
novel for $q<\infty $, and gives for the first time a Fourier-multiplier regularity
estimate. When $q=\infty$ and $\frac{1}{r} \ne \frac{1}{p}$, it is known that logarithmic Besov regularity
\cite{ben2019convergence,brue2021sharp,brue2021advection,meyer2022propagation}
extends to a fractional power $\frac{p}{r}$ of the logarithm, with a
polynomial in time bound. One might expect a similar generalisation of our
result, when $\frac{1}{r} \ne \frac{1}{p} + \frac{1}{q}$. However, for non-Besov estimates this remains an open
problem except in the case $r=2, q=\infty$ when the Besov and Fourier-multiplier estimate agree.
        

    \setcounter{subsection}{3}
    \subsection{Vanishing diffusion and convergence}\label{sec:convergence}
        As well as mixing and propagation of regularity, the tools in this paper
can also be used to quantify convergence rates for vanishing diffusion
and mollification limits. Such problems have been studied in the literature,
but either require additional smoothness of the initial datum
\cite{brue2021advection,meyer2022propagation}, or rely on optimal transport
to quantify convergence rates
\cite{seis2017quantitative,seis2018optimal}. In contrast, we prove
logarithmic convergence rates
for all initial data in the DiPerna-Lions
class $\rho _{0}\in L_{x}^{q}$, by convergence of frequency
cutoffs.


        Suppose $\rho_\nu(x,t):\mathbb{R}^d \times[0,T]\to\mathbb{R}$ with $\rho_\nu(x,t) \in L_t^\infty L_x^q$ solves the advection-diffusion equation with the same initial datum $\rho_0(x)\in L_x^q$,
        \begin{align*}
            \begin{cases}
                \frac{\partial \rho_\nu}{\partial t} + \nabla\cdot(u\rho_\nu) - \nu\Delta\rho_\nu = 0, \\
                \rho_\nu(x,0)=\rho_0(x).
            \end{cases}
        \end{align*}
        Then, for suitably small diffusion $\nu > 0$, we have the bound
        \begin{equation}\label{eq:stability1intro}
            \left\|(\rho-\rho_\nu)*\varphi_\delta\right\|_{L_t^\infty L_x^r} \le C_{p,q} \left\|\rho_0\right\|_{L_x^q}\left\|\nabla u\right\|_{L_t^1 L_x^p}\left(\log\left(\frac{\delta^2}{\nu T}\left(\frac{\left\|\rho_0\right\|_{L_x^q}}{\left\|\rho_0\right\|_{L_x^r}}\left\|\nabla u\right\|_{L_t^1 L_x^p}\right)\right)\right)^{-1},
        \end{equation}
        where $\varphi_\delta(x)$ denotes the frequency cutoff \eqref{eq:frequencycutoffintro}. See Proposition \ref{thm:vanishingviscosity} for the exact statement of the theorem. See Corollary \ref{cor:vanishingviscosity} for convergence of the passive scalar in the space $W_x^{-1,r}$.
        
        Similarly, let $\bar{\rho}(x,t):\mathbb{R}^d\times[0,T]\to\mathbb{R}$ with $\bar{\rho}(x,t) \in L_t^\infty L_x^q$ solve the transport equation
        \begin{gather*}
            \begin{cases}
                \frac{\partial \bar{\rho}}{\partial t} + \nabla\cdot(\bar{u}\bar{\rho}) = 0, \\
                \bar{\rho}(x,0)=\rho_0(x),
            \end{cases}
        \end{gather*}
        for the same initial datum $\rho_0(x)\in L_x^q$, along a measureable vector field $\bar{u}(x,t)\in L_t^1 L_x^p$. Then, if $\left\|u-\bar{u}\right\|_{L_t^1L_x^p}\le\delta\left\|\nabla u\right\|_{L_t^1 L_x^p}$,
        \begin{equation}\label{eq:stability2intro}
            \left\|(\rho-\bar{\rho})*\varphi_\delta\right\|_{L_t^\infty L_x^r} \le C_{p,q}\left\|\bar{\rho}\right\|_{L_t^\infty L_x^q}\left\|\nabla u\right\|_{L_t^1L_x^p}\left(\log\left(\frac{\delta\left\|\nabla u\right\|_{L_t^1 L_x^p}}{\left\|u-\bar{u}\right\|_{L_t^1L_x^p}}\right)\right)^{-1},
        \end{equation}
        where $\varphi_\delta(x)$ again denotes the frequency cutoff \eqref{eq:frequencycutoffintro}. See Proposition \ref{thm:stability} for the exact statement of the theorem. See Corollary \ref{cor:stability} for convergence of the passive scalar in the space $W_x^{-1,r}$. 

        In \cite{brue2021advection,meyer2022propagation} the authors show logarithmic convergence rates similar to \eqref{eq:stability1intro}, \eqref{eq:stability2intro} in the strong topology $L_x^2$ for more regular initial data. The main improvement of Corollary \ref{cor:vanishingviscosity}, \ref{cor:stability} is that we obtain convergence rates for any initial data. One can combine \eqref{eq:stability1intro}, \eqref{eq:stability2intro} with Proposition \ref{thm:propagation} to similarly recover strong convergence for more regular initial data, or by estimating the remainder $\|(\rho-\rho_\nu) -(\rho-\rho_\nu) * \varphi_\delta\|_{L_t^\infty L_x^r}$ using Theorem \ref{thm:quantitative2}. We note, however, that the convergence rates themselves do not improve by measuring them in weak norms compared to \cite{brue2021advection,meyer2022propagation}.



    \setcounter{subsection}{4}
    \subsection{Connection with the standard commutator}\label{sec:diperna}
        Finally, we discuss how the key harmonic estimates in this paper are related to quantitative estimates on the decay of the standard DiPerna-Lions commutator, as is necessary for the well-posedness theory \cite{diperna1989ordinary} for transport by Sobolev vector fields. We note that the DiPerna-Lions commutator arises as a special case of the commutator in Lemma \ref{thm:calderon} by taking the kernel $K(x)$ to be a mollifier $K(x)=\varphi_\delta(x)$. The resulting commutator $R_\delta(x,t)$ controls the the transport of $(\rho*\varphi_\delta)(x,t)$, and the key observation of \cite{diperna1989ordinary} is that $\left\|R_\delta(\cdot,t)\right\|_{L_x^r}$ vanishes as $\delta\to0$, from which well-posedness follows (for all $p\ge1$). Quantitative estimates on the decay of the DiPerna-Lions commutator should, therefore, give quantitative estimates for the well-posedness of the Cauchy problem, as in Theorem \ref{thm:quantitative2}. This was the original motivation for this work. Indeed, the key harmonic estimate in the proof of Theorem \ref{thm:quantitative2} is equivalent to an \textit{integral}, or \textit{averaged}, decay rate of the DiPerna-Lions commutator. In particular, we show that
        \begin{equation}\label{eq:dipernalionsintegral}
            \left\|\int_0^\infty R_\delta(x,t)\;\frac{d\delta}{\delta} \right\|_{L_x^r} \le C_{p,q}\left\|\rho(\cdot,t)\right\|_{L_x^q}\left\|\nabla u(\cdot,t)\right\|_{L_x^p},
        \end{equation}
        see Proposition \ref{thm:expdecay} for details. When sharing a version of this manuscript, we were made aware of an unpublished result similar to \eqref{eq:dipernalionsintegral} by Brué, Colombo and De Philippis \cite{eliaprivate} for the special case of the Gaussian commutator. Since the integral against $\frac{d\delta}{\delta}$ is logarithmically singular, \eqref{eq:dipernalionsintegral} implies some averaged control on the pointwise decay of $R_\delta(x,t)\to0$ of the order $|\log \delta|^{-1}$. While this is not decay of the norm $\left\|R_\delta(\cdot,t)\right\|_{L_x^r}\to0$ as in \cite{diperna1989ordinary}, it is sufficient to prove quantitative estimates, and in particular our main result Theorem \ref{thm:quantitative2}. We now ask if one can improve this estimate \eqref{eq:dipernalionsintegral} to a more standard decay rate for the DiPerna-Lions commutator, potentially leading to an improvement of Theorem \ref{thm:quantitative2}.

        By different methods, we show that one has the bound
        \begin{equation*}
            \left(\int_0^\infty \left\|R_\delta(\cdot,t)\right\|_{L_x^r}^{\max(p,q)}\;\frac{d\delta}{\delta}\right)^\frac{1}{\max(p,q)} \le C_{p,q}\left\|\rho(\cdot,t)\right\|_{L_x^q}\left\|\nabla u(\cdot,t)\right\|_{L_x^p},
        \end{equation*}
        see Proposition \ref{thm:besovdecay}, and Corollary \ref{thm:lpdecay} for details. Rather than a quantitative decay rate, which does not exist (see point 1 of Proposition \ref{thm:counterexamples}), this gives an integral decay rate for the \textit{norm} $\left\|R_\delta(\cdot,t)\right\|_{L_x^r}\to0$, at the expense of the \textit{rate} $\max(p,q)$. Curiously, this rate is sharp when $p,q=2$, see point 2 of Proposition \ref{thm:counterexamples}. As a consequence, this approach leads to suboptimal stability and mixing rates compared to the estimate \eqref{eq:dipernalionsintegral}, and so we do not expand on it further. However, we include it for the curious reader and also because in the follow-up work \cite{BVpaper} we extend this approach to give quantitative estimates for vector fields in $W_x^{1,p}$ with $p=1$. Indeed, it was an open question to extend any of the quantitative results in the literature to vector fields in $W_x^{1,1}$ or $BV_x$.

    \subsection{Outline of the paper}
        In Section \ref{section2}, we give the proofs of the key harmonic estimate, Lemma \ref{thm:calderon}, and our main result on the stability of the Cauchy problem, Theorem \ref{thm:quantitative2}.
        
        In Section \ref{section3}, we use Theorem \ref{thm:quantitative2} and Lemma \ref{thm:calderon} to prove the quantitative estimates on mixing (Theorem \ref{thm:mixing}), propagation of regularity (Proposition \ref{thm:propagation}), convergence of vanishing diffusion (Proposition \ref{thm:vanishingviscosity}, Corollary \ref{cor:vanishingviscosity}) and mollification of the vector field (Proposition \ref{thm:stability}, Corollary \ref{cor:stability}).
        
        Finally, in Section \ref{section4}, we discuss connections with the standard DiPerna-Lions commutator, proving in particular the integral point-wise decay rate Proposition \ref{thm:expdecay}, and the integral norm decay rate Proposition \ref{thm:besovdecay}.

\section{Quantitative well-posedness estimates}\label{section2}
    For any convolution kernel $K(x):\mathbb{R}^d\to\mathbb{R}$, consider the commutator $R(x,t)$ which gives the error in the transport of $(\rho*K)(x,t)$:
    \begin{equation*}
        \left(\frac{\partial}{\partial t} + u(x,t)\cdot\nabla\right)(\rho * K)(x,t) = R(x,t),
    \end{equation*}
    i.e.
    \begin{equation*}
        R(x,t) = u(x,t)\cdot(\rho*\nabla K)(x,t) - \nabla\cdot((u\rho)*K)(x,t).
    \end{equation*}
    By the Chris-Journ\'e singular integral estimate in \cite{seeger2019multilinear}, this commutator admits the following bound:
    \begin{lemma}\label{thm:calderon}
        Consider a scalar $\rho(x):\mathbb{R}^d\to\mathbb{R}$ with $\rho(x)\in L_x^q$, and a divergence-free vector field $u(x):\mathbb{R}^d\to\mathbb{R}^d$, with $\nabla u(x)\in L_x^p$. Let $\frac{1}{r}=\frac{1}{p}+\frac{1}{q}$, with $1<p,q\le\infty$, $1\le r<\infty$.
        
        For any convolution kernel $K(x):\mathbb{R}^d\to\mathbb{R}$, define the commutator
        \begin{equation*}
            R(x) = u\cdot\nabla(\rho*K) - \nabla\cdot\left((u\rho)*K\right),
        \end{equation*}
        then
        \begin{equation*}
            \left\|R\right\|_{L_x^r}\le C_{p,q} \left\|\rho\right\|_{L_x^q}\left\|\nabla u\right\|_{L_x^p}\left(\sum_{i,j=1}^d\left\|\frac{\partial(x_jK)}{\partial x_i}\right\|_{CZ}\right),
        \end{equation*}
        where the constant $C_{p,q}>0$ depends only on the parameters $1<p,q\le\infty$ and the dimension $d\ge1$, and where $\left\|\cdot\right\|_{CZ}$ refers to the Calderon-Zygmund norm, i.e. the minimum constant $\|K\|_{CZ}$ such that
        \begin{equation*}
            \|K\|_{CZ} = \mathrm{min}\left\{\left\||x|^dK\right\|_{L_x^\infty},\left\||x|^{d+1}\nabla K\right\|_{L_x^\infty},\left\|\hat{K}\right\|_{L_\xi^\infty}\right\},
        \end{equation*}
        where $\hat{K}(\xi):\mathbb{R}^d\to\mathbb{C}$ denotes the Fourier transform.
    \end{lemma}
    \begin{proof}
        We have, as for the DiPerna-Lions commutator \cite{diperna1989ordinary},
        \begin{align*}
            R(x) & = \int_{\mathbb{R}^d}\rho(y)(u(x)-u(y))\cdot\nabla K(x-y)\;dy \\
            & = \int_{\mathbb{R}^d} \rho(x-h)\sum^d_{i,j=1}\left(\int_0^1h_j\frac{\partial u_i}{\partial x_j}(x-th)\;dt\right) \frac{\partial K}{\partial h_i}(h)\;dh \\
            & = \begin{aligned}[t]
                & \sum^d_{i,j=1}\int_{\mathbb{R}^d} \rho(x-h)\left(\int_0^1\frac{\partial u_i}{\partial x_j}(x-th)\;dt\right) \frac{\partial (h_j K)}{\partial h_i}(h)\;dh \\
                & - \int_{\mathbb{R}^d} \rho(x-h)\left(\int_0^1(\nabla\cdot u)(x-th)\;dt\right) K(h)\;dh.
            \end{aligned}
        \end{align*}

       hence since $\nabla\cdot u(x)=0$, we have
        \begin{equation*}
            \left\|R\right\|_{L_x^r}\le \sum_{i,j=1}^d\left\|\int_{\mathbb{R}^d} \rho(x-h)\left(\int_0^1\frac{\partial u_i}{\partial x_j}(x-th)\;dt\right) \frac{\partial (h_j K)}{\partial h_i}(h)\;dh\right\|_{L_x^r},
        \end{equation*}
        and so we conclude by the key harmonic analysis estimate \cite{seeger2019multilinear}*{Theorem 1.1}.
    \end{proof}

    From this, we deduce the main result of this paper, which unifies various stability, mixing, and well-posedness estimates:

    \begin{definition}[Mollifier]\label{def:mollifier}
        We say $\varphi(x)\in L^1(\mathbb{R}^d;\mathbb{R})$ is a standard mollifier if
        \begin{equation*}
            \int_{\mathbb{R}^d}\varphi(x)\;dx = 1 \text{ and define }\varphi_\delta(x) \in L^1(\mathbb{R}^d), \ \varphi_\delta(x) = \delta^{-d}\varphi\left(\frac{x}{\delta}\right),
        \end{equation*}
        for all $\delta > 0$. Note that one often asks for additional smoothness of $\varphi(x)$ such as $\varphi(x) \in C_c^\infty(\mathbb{R}^d)$. We will make this explicit when it is required.
    \end{definition}

    \begin{definition}[Frequency cutoff]\label{def:frequencycutoff}
        We say a mollifier $\varphi(x)\in \mathcal{S}(\mathbb{R}^d;\mathbb{R})$ (a Schwartz function) is a frequency-cutoff if its Fourier transform $\hat{\varphi}(\xi)\in \mathcal{S}(\mathbb{R}^d;\mathbb{C})$ satisfies the following:
        \begin{equation}\label{eq:frequencycutoff}
            \begin{gathered}
                \hat{\varphi}(\xi)=1 \text{ for } |\xi|\le1, \\
                \hat{\varphi}(\xi)=0 \text{ for } |\xi|\ge 2.
            \end{gathered}
        \end{equation}
    \end{definition}

    We now give the main result of this paper, a quantitative well-posedness result for passive scalar transport, which quantifies the propagation of the smallness of the initial data. Continuous dependence on the initial data in the strong topology is standard by the theory of renormalised solutions, \cite{diperna1989ordinary}. In contrast, using Lemma \ref{thm:calderon}, we show continuous dependence on the initial data in the weak topology, quantified by frequency cutoffs:

    \begin{theorem}\label{thm:quantitative2}
        Consider a weak solution $\rho(x,t):\mathbb{R}^d\times[0,T]\to\mathbb{R}$, with $\rho(x,t)\in L_t^\infty L_x^q$, and initial datum $\rho_0(x)\in L_x^q$, to the forced transport equation
        \begin{equation*}
            \begin{cases}
                \frac{\partial \rho}{\partial t}(x,t)+\nabla \cdot(u(x,t)\rho(x,t))=f(x,t), \\
                \rho(x,0)=\rho_0(x),
            \end{cases}
        \end{equation*}
        along a divergence-free vector field $u(x,t):\mathbb{R}^d\times[0,T]\to\mathbb{R}^d$, with $\nabla u(x,t)\in L_t^1L_x^p$, for $\frac{1}{p}+\frac{1}{q}=\frac{1}{r}$, with $1<p,q\le\infty$, $1\le r<\infty$, and distributional force $f(x,t)\in\mathcal{S}'(\mathbb{R}^d\times\mathbb{R};\mathbb{R})$ a Schwartz distribution supported on $\mathbb{R}^d\times[0,T]$.

        Fix a mollifier/frequency cutoff $\varphi(x)\in\mathcal{S}(\mathbb{R}^d;\mathbb{R})$ as in Definition \ref{def:frequencycutoff}.

        Then for any $\delta>0$ and $\kappa > 0$ we have the following quantitative bound:
        
        If
        \begin{equation*}
            \delta_0 \le \frac{1}{4}\delta \exp\left(-\frac{C_{p,q}}{\kappa}\left\|\rho\right\|_{L_t^\infty L_x^q}\left\|\nabla u\right\|_{L_t^1 L_x^p}\right),
        \end{equation*}
        then
        \begin{equation*}
            \left\|\rho*\varphi_{\delta}\right\|_{L_t^\infty L_x^r} \le \kappa + A\left\|\rho_0*\varphi_{\delta_0}\right\|_{L_x^r} + A\left\|f*\varphi_{\delta_0}\right\|_{L_t^1 L_x^r},
        \end{equation*}
        where the constant $C_{p,q} > 0$ depends only on the parameters $1<p,q\le\infty$, the choice of mollifier $\varphi(x)\in\mathcal{S}(\mathbb{R}^d;\mathbb{R})$, and the dimension $d\ge1$, while $A>0$ depends only on the choice of mollifier $\varphi(x)\in\mathcal{S}(\mathbb{R}^d;\mathbb{R})$ and the dimension $d\ge1$.
    \end{theorem}
    \begin{proof}
        Throughout, $\lesssim$ will denote less than or equal to up to a constant depending only on the choice of mollifier $\varphi(x)\in\mathcal{S}(\mathbb{R}^d;\mathbb{R})$, and the dimension $d\ge1$.

        We consider the following convolution kernel in Lemma \ref{thm:calderon}:
        \begin{equation}\label{eq:kernel}
            K(x)=\int_{\delta_1}^{\delta_2}\varphi_\delta(x)\;\frac{d\delta}{\delta} \in \mathcal{S}(\mathbb{R}^d;\mathbb{R}),
        \end{equation}
        for some $0<\delta_1<\delta_2$ to be chosen later.

        Let $\psi(x,t)\in C_c^\infty(\mathbb{R}^d\times[0,T];\mathbb{R})$ be any test function. Let $\phi(x,t)\in L_t^\infty L_x^p$ be a weak solution to the \textit{backwards} transport equation with force $\psi(x,t)$:
        \begin{equation}\label{eq:testtransport}
            \begin{gathered}
                \frac{\partial\phi}{\partial t}(x,t)+\nabla\cdot(u(x,t)\phi(x,t))=\psi(x,t), \\
                \phi(x,T)=0.
            \end{gathered}
        \end{equation}
        satisfying the standard energy estimate
        \begin{equation}\label{eq:testenergyestimate}
            \left\|\phi\right\|_{L_t^\infty L_x^{r'}}\le\left\|\psi\right\|_{L_t^1L_x^{r'}},
        \end{equation}
        for $\frac{1}{r'}+\frac{1}{r}=1$, see for instance the lecture notes \cite{ambrosio2014continuity}.

        Convolving this equation with $\bar{K}(x) = K(-x)$ gives
        \begin{equation}\label{eq:testtransportmollified}
            \frac{\partial}{\partial t}(\phi*\bar{K})(x,t)+\nabla\cdot(u(x,t)(\phi*\bar{K})(x,t))=R(x,t) + (\psi*\bar{K})(x,t),
        \end{equation}
        in terms of the commutator
        \begin{align*}
            R(x,t) & = u(x,t)\cdot\nabla(\phi*\bar{K})(x,t)-(\nabla\cdot(u\phi)*\bar{K})(x,t) \\
            & = \int_{\mathbb{R}^d}\phi(y,t)(u(x,t)-u(y,t))\cdot\nabla \bar{K}(x-y)\;dy,
        \end{align*}
        
        Applying Lemma \ref{thm:calderon}, for $\frac{1}{q'}+\frac{1}{q}=1$, $\frac{1}{r'}+\frac{1}{r}=1$, we have the bound
        \begin{align}\label{eq:remainderbound0}
            \left\|R(x,t)\right\|_{L_t^1L_x^{q'}} \le C_{p,q} \left\|\phi\right\|_{L_t^\infty L_x^{r'}}\left\|\nabla u\right\|_{L_t^1L_x^p}\left(\sum_{i,j=1}^d\left\|\frac{\partial(x_j\bar{K})}{\partial x_i}\right\|_{CZ}\right),
        \end{align}
        for some constant $C_{p,q}>0$ depending only on the parameters $1<p,q\le\infty$ and the dimension $d\ge1$.

        To estimate this Calderon-Zygmund norm, we note that
        \begin{equation*}
            \left\|\frac{\partial(x_j\bar{K})}{\partial x_i}\right\|_{CZ} = \left\|\frac{\partial(x_jK)}{\partial x_i}\right\|_{CZ},
        \end{equation*}
        and so define
        \begin{align*}
            K_{i,j}'(x) & = \frac{\partial (x_j K)}{\partial x_i}(x) = \int_{\delta_1}^{\delta_2}\frac{\partial (x_j\varphi)}{\partial x_i}\left(\frac{x}{\delta}\right)\;\frac{d\delta}{\delta^{d+1}}.
        \end{align*}
        
        We then have the following estimates which imply that $K_{i,j}'$ is a Calderon-Zygmund kernel:
        {\allowdisplaybreaks
        \begin{gather}
            \begin{aligned}[t]
                |K_{i,j}'(x)| & \le \left\|(1+|x|^{d+1})\frac{\partial (x_j\varphi)}{\partial x_i}\right\|_{L_x^\infty}\left(\int_{|x|}^\infty \delta^{-d-1} d\delta + \int_0^{|x|} |x|^{-d-1} d\delta \right) \\
                & \lesssim |x|^{-d},
            \end{aligned} \label{eq:calderonestimate1} \\
            \begin{aligned}[t]
                \left|\frac{\partial K_{i,j}'}{\partial x_k}(x)\right| & = \left|\int_{\delta_1}^{\delta_2} \frac{\partial^2 (x_j\varphi)}{\partial x_k \partial x_i}\left(\frac{x}{\delta}\right)\;\frac{d\delta}{\delta^{d+2}}\right| \\
                & \le \left\|(1+|x|^{d+2})\frac{\partial^2 (x_j\varphi)}{\partial x_k \partial x_i}\right\|_{L_x^\infty}\left(\int_{|x|}^\infty \delta^{-d-2} d\delta + \int_0^{|x|} |x|^{-d-2} d\delta \right) \\
                & \lesssim |x|^{-d-1},
            \end{aligned} \label{eq:calderonestimate2} \\
            \begin{aligned}[t]
                |\hat{K}_{i,j}'(\xi)| & = \left| \int_{\mathbb{R}^d} e^{i\xi\cdot x} \int_{\delta_1}^{\delta_2} \frac{\partial (x_j\varphi)}{\partial x_i}\left(\frac{x}{\delta}\right) \; \frac{d\delta}{\delta^{d+1}} \; dx \right| \\
                & = \left|\int_{\delta_1}^{\delta_2} \left(\int_{\mathbb{R}^d} e^{i \delta \xi\cdot h} \frac{\partial (x_j\varphi)}{\partial x_i}(h)\; dh \right) \frac{d\delta}{\delta} \right| \\
                & = \left|\int_{\delta_1}^{\delta_2} -i\delta \xi_i \frac{\partial\hat{\varphi}}{\partial \xi_j}(\delta \xi) \; \frac{d\delta}{\delta} \right| \\
                & \le \int_0^\infty |\xi_i| \left|\frac{\partial\hat{\varphi}}{\partial \xi_j}\right|(\delta \xi)\; d\delta \\
                & \le \left\|(1+|\xi_i|^2)\frac{\partial\hat{\varphi}}{\partial \xi_j}\right\|_{L_x^\infty}\left(\int_0^{|\xi_i|^{-1}} |\xi_i| \; d\delta + \int_{|\xi_i|^{-1}}^\infty |\xi_i|^{-1} \delta^{-2}\; d\delta\right) \\
                & \lesssim 1.
            \end{aligned} \label{eq:calderonestimate3}
        \end{gather}
        }
        where, in particular, none of the bounds depends on $0<\delta_1<\delta_2$. That is
        \begin{equation*}
            \left\|K_{i,j}'\right\|_{CZ} \lesssim 1.
        \end{equation*}

        Therefore, by Lemma \ref{thm:calderon} and the energy estimate $\left\|\phi\right\|_{L_t^\infty L_x^{r'}}\le\left\|\psi\right\|_{L_t^1L_x^{r'}}$ in \eqref{eq:testenergyestimate}, then \eqref{eq:remainderbound0} becomes
        \begin{equation}\label{eq:remainderbound}
            \left\|R(x,t)\right\|_{L_t^1L_x^{q'}} \lesssim C_{p,q} \left\|\psi\right\|_{L_t^1 L_x^{r'}}\left\|\nabla u\right\|_{L_t^1L_x^p}.
        \end{equation}

        By the expression for the kernel $K(x)$ in \eqref{eq:kernel}, we estimate
        \begin{align}
            \left\|\phi*\bar{K}\right\|_{L_t^\infty L_x^{r'}} & \le \left\|K\right\|_{L_x^1} \left\|\phi\right\|_{L_t^\infty L_x^{r'}} \nonumber \\
            & \lesssim \log\left(\frac{\delta_2}{\delta_1}\right)\left\|\psi\right\|_{L_t^1L_x^{r'}} \label{eq:weakcontinuitybound}
        \end{align}

        Recall that $\rho(x,t)$ is a weak solution to the following forced transport equation
        \begin{equation}\label{eq:forcedtransport}
            \frac{\partial \rho}{\partial t}(x,t)+\nabla \cdot(u(x,t)\rho(x,t))=f(x,t),
        \end{equation}
        with distributional force $f(x,t)\in\mathcal{S}'(\mathbb{R}^d\times\mathbb{R};\mathbb{R})$ supported on $\mathbb{R}^d\times[0,T]$.
        
        We shall eventually take $\delta_1=2\delta_0$. Assume that the distribution $(f*\varphi_{2^{-1}\delta_1})(x,t)$ has finite $L_t^1 L_x^r$-norm, else the statement of the theorem is trivial. Then also $(f*\varphi_{2^{-1}\delta_1}*K)(x,t)\in L_t^1 L_x^r$. 

        Since the mollifier $\varphi(x)\in\mathcal{S}(\mathbb{R}^d;\mathbb{R})$ is a frequency-cutoff \eqref{eq:frequencycutoff}, we see that for all $\delta\ge2\delta_2$, and $\delta_1\le\delta'\le\delta_2$
        \begin{equation}\label{eq:fourierproperty}
            \begin{gathered}
                \varphi_{\delta} = \varphi_{\delta'} * \varphi_{\delta}, \\
                \varphi_{\delta'} = \varphi_{2^{-1}\delta_1} * \varphi_{\delta'}.
            \end{gathered}
        \end{equation}
        In particular, by integrating over $\int_{\delta_1}^{\delta_2}\frac{d\delta'}{\delta'}$,
        \begin{equation}\label{eq:fourierproperty2}
            \begin{gathered}
                \log\left(\frac{\delta_2}{\delta_1}\right) \varphi_{\delta}(x) = (K*\varphi_{\delta})(x), \\
                K(x) = (\varphi_{2^{-1}\delta_1}*K)(x).
            \end{gathered}
        \end{equation}

        Since $(f*\varphi_{2^{-1}\delta_1}*K)(x,t)\in L_t^1 L_x^r$, then $(f*K)(x,t)\in L_t^1 L_x^r$.
        
        Therefore, we may take $(\phi*\bar{K})(x,t)$ as a test function to \eqref{eq:forcedtransport}:
        \begin{align*}
            & \int_{\mathbb{R}^d\times[0,T]}\rho(x,t)\left(\frac{\partial}{\partial t}(\phi*\bar{K})(x,t)+u(x,t)\cdot\nabla(\phi*\bar{K})(x,t)\right)\;dxdt \\
            & \qquad = - \int_{\mathbb{R}^d}\rho_0(x)(\phi*\bar{K})(x,0)\;dx - \int_{\mathbb{R}^d;\mathbb{R}}(f*K)(x,t)\phi(x,t)\;dxdt.
        \end{align*}

        and so, using \eqref{eq:testtransportmollified},
        \begin{align*}
            & \int_{\mathbb{R}^d\times[0,T]}\rho(x,t)(\psi*\bar{K})(x,t)\;dxdt + \int_{\mathbb{R}^d\times[0,T]}\rho(x,t)R(x,t)\;dxdt \\
            & \qquad = - \int_{\mathbb{R}^d}\rho_0(x)(\phi*\bar{K})(x,0)\;dx - \int_{\mathbb{R}^d\times[0,T]}(f*K)(x,t)\phi(x,t)\;dxdt.
        \end{align*}

        By \eqref{eq:fourierproperty2}, we also have $\bar{K}(x)=(\bar{K}*\varphi_{2^{-1}\delta_1})(x)$, and so
        \begin{align*}
            & \int_{\mathbb{R}^d\times[0,T]}(\rho*K)(x,t)\psi(x,t)\;dxdt + \int_{\mathbb{R}^d\times[0,T]}\rho(x,t)R(x,t)\;dxdt \\
            & \qquad = - \int_{\mathbb{R}^d}(\rho_0*\varphi_{2^{-1}\delta_1})(x)(\phi*\bar{K})(x,0)\;dx - \int_{\mathbb{R}^d;\mathbb{R}}(f*\varphi_{2^{-1}\delta_1})(x,t)(\phi*\bar{K})(x,t)\;dxdt.
        \end{align*}
        
        Therefore, by \eqref{eq:remainderbound}, \eqref{eq:weakcontinuitybound} we have the bound
        \begin{align*}
            & \int_{\mathbb{R}^d\times[0,T]}(\rho*K)(x,t)\psi(x,t)\;dxdt \\
            & \qquad \le \begin{aligned}[t]
                & \left\|\rho\right\|_{L_t^\infty L_x^q}\left\|R\right\|_{L_t^1 L_x^{q'}} + \left\|\rho_0*\varphi_{2^{-1}\delta_1}\right\|_{L_x^r}\left\|K\right\|_{L_x^1}\left\|\phi\right\|_{L_t^\infty L_x^{r'}} \\
                & + \left\|f*\varphi_{2^{-1}\delta_1}\right\|_{L_t^1 L_x^r}\left\|\phi*\bar{K}\right\|_{L_t^\infty L_x^{r'}}
            \end{aligned} \\
            & \qquad \lesssim \begin{aligned}[t]
                & C_{p,q}\left\|\rho\right\|_{L_t^\infty L_x^q}\left\|\psi\right\|_{L_t^1 L_x^{r'}}\left\|\nabla u\right\|_{L_t^1 L_x^p} + \log\left(\frac{\delta_2}{\delta_1}\right)\left\|\rho_0*\varphi_{2^{-1}\delta_1}\right\|_{L_x^r}\left\|\psi\right\|_{L_t^1 L_x^{r'}} \\
                & + \log\left(\frac{\delta_2}{\delta_1}\right)\left\|f*\varphi_{2^{-1}\delta_1}\right\|_{L_t^1 L_x^r}\left\|\psi\right\|_{L_t^1L_x^{r'}},
            \end{aligned}
        \end{align*}
        and since $\psi(x,t)\in C_c^\infty(\mathbb{R}^d\times[0,T];\mathbb{R})$ is arbitrary, this becomes
        \begin{equation*}
            \left\|\rho*K\right\|_{L_t^\infty L_x^r} \lesssim C_{p,q}\left\|\rho\right\|_{L_t^\infty L_x^q}\left\|\nabla u\right\|_{L_t^1 L_x^p} + \log\left(\frac{\delta_2}{\delta_1}\right)\left\|\rho_0*\varphi_{2^{-1}\delta_1}\right\|_{L_x^r} + \log\left(\frac{\delta_2}{\delta_1}\right)\left\|f*\varphi_{2^{-1}\delta_1}\right\|_{L_t^1 L_x^r}.
        \end{equation*}

        Then by \eqref{eq:fourierproperty2}, for any $\delta\ge2\delta_2$,
        \begin{align*}
            \left\|\rho*\varphi_{\delta}\right\|_{L_t^\infty L_x^r} & \le \frac{1}{\log\left(\frac{\delta_2}{\delta_1}\right)}\left\|\rho*K\right\|_{L_t^\infty L_x^r} \\
            & \lesssim \frac{C_{p,q}}{\log\left(\frac{\delta_2}{\delta_1}\right)}\left\|\rho\right\|_{L_t^\infty L_x^q}\left\|\nabla u\right\|_{L_t^1 L_x^p} + \left\|\rho_0*\varphi_{2^{-1}\delta_1}\right\|_{L_x^r} + \left\|f*\varphi_{2^{-1}\delta_1}\right\|_{L_t^1 L_x^r}.
        \end{align*}

        We now pick $0<\delta_1<\delta_2$. For any $\kappa>0$, let
        \begin{gather*}
            \delta_2=\delta_1\exp\left(\frac{C_{p,q}}{\kappa}\left\|\rho\right\|_{L_t^\infty L_x^q}\left\|\nabla u\right\|_{L_t^1 L_x^p}\right).
        \end{gather*}

        Then for any
        \begin{equation*}
            \delta \ge 4(2^{-1}\delta_1) \exp\left(\frac{C_{p,q}}{\kappa}\left\|\rho\right\|_{L_t^\infty L_x^q}\left\|\nabla u\right\|_{L_t^1 L_x^p}\right),
        \end{equation*}
        we have the following bounds:
        \begin{equation*}
            \left\|\rho*\varphi_{\delta}\right\|_{L_t^\infty L_x^r} \lesssim \kappa + \left\|\rho_0*\varphi_{2^{-1}\delta_1}\right\|_{L_x^r} + \left\|f*\varphi_{2^{-1}\delta_1}\right\|_{L_t^1 L_x^r}.
        \end{equation*}
        The result now follows by letting $\delta_1=2\delta_0$.
    \end{proof}

\section{Mixing, regularity, and convergence}\label{section3}
    We now apply Lemma \ref{thm:calderon} and Theorem \ref{thm:quantitative2} to four well-studied problems in the theory of passive scalar transport: mixing, approximation of the vector field, the vanishing diffusion limit, and regularity of the passive scalar.

    \begin{definition}[$W_x^{-1,r}$]\label{def:scaledsobolev}
        For $\frac{1}{r'}+\frac{1}{r}=1$, we define the norm $W_x^{-1,r}$ as the dual to
        \begin{equation*}
            \left\|\phi\right\|_{W_x^{1,r'}}=\left\|\phi\right\|_{L_x^{r'}} + \left\|\nabla \phi\right\|_{L_x^{r'}}.
        \end{equation*}
    \end{definition}

    We now apply Theorem \ref{thm:quantitative2} to give an exponential lower bound on the decay of the negative Sobolev norm $\left\|\rho(\cdot,t)\right\|_{W_x^{-1,r}}$, also studied in \cite{seis2013maximal,iyer2014lower,leger2018new}. The norm $W_x^{-1,r}$ measures the average length scale of the passive scalar at time $t>0$, and so quantifies how mixed is the passive scalar:

    \begin{theorem}[Mixing bound]\label{thm:mixing}
        Consider a weak solution $\rho(x,t):\mathbb{R}^d\times[0,T]\to\mathbb{R}$, with $\rho(x,t)\in L_t^\infty L_x^q$, and initial datum $\rho_0(x)\in L_x^q$, to the transport equation
        \begin{equation*}
            \begin{cases}
                \frac{\partial \rho}{\partial t}(x,t)+\nabla \cdot(u(x,t)\rho(x,t))=0, \\
                \rho(x,0)=\rho_0(x),
            \end{cases}
        \end{equation*}
        along a divergence-free $u(x,t):\mathbb{R}^d\times[0,T]\to\mathbb{R}^d$, with $\nabla u(x,t)\in L_t^1L_x^p$, for $\frac{1}{p}+\frac{1}{q}=\frac{1}{r}$, with $1<p,q\le\infty$, $1\le r<\infty$.

        Let $\varphi(x)\in L^1(\mathbb{R}^d;\mathbb{R})$ be any \textit{positive} mollifier. For any $0<\delta \le 1$, such that
        
        \begin{equation*}
            \left\|\rho_0*\varphi_\delta\right\|_{L_x^r}\ge\frac{1}{2}\left\|\rho_0\right\|_{L_x^r},
        \end{equation*}
        then
        \begin{equation*}
            \left\|\rho(\cdot,T)\right\|_{W_x^{-1,r}}\ge \delta A\left\|\rho_0\right\|_{L_x^r}\exp\left(-C_{p,q}\frac{\left\|\rho_0\right\|_{L_x^q}}{\left\|\rho_0\right\|_{L_x^r}}\left\|\nabla u\right\|_{L_t^1L_x^p}\right),
        \end{equation*}
        where the constant $C_{p,q}>0$ depends only on the parameters $1<p,q\le\infty$ and the dimension $d\ge1$, while $A>0$ depends on the choice of positive mollifier $\varphi(x)\in L^1(\mathbb{R}^d;\mathbb{R})$ and the dimension $d\ge1$.
    \end{theorem}
    \begin{proof}
        For such a weak solution $\rho(x,t)$, define the time reverse $\tilde{\rho}(x,t)=\rho(x,T-t)$, and similar $\tilde{u}(x,t)=-u(x,T-t)$. Then $\tilde{\rho}(x,t)$ is a weak solution to the transport equation along $\tilde{u}(x,t)$ with initial datum $\tilde{\rho}_0(x)=\rho(x,T)$ satisfying $\left\|\tilde{\rho}_0\right\|_{L_x^q}=\left\|\rho_0\right\|_{L_x^q}$, $\left\|\tilde{\rho}_0\right\|_{L_x^r}=\left\|\rho_0\right\|_{L_x^r}$ by the renormalization property, see for instance \cite{diperna1989ordinary}.
        
        Therefore, by Theorem \ref{thm:quantitative2}, for all $\delta_0>0$ and any $\kappa > 0$, we have the following quantitative bound on the propagation of length scales. For
        \begin{align*}
            \delta' & = 4\delta_0 \exp\left(\frac{C_{p,q}}{\kappa}\left\|\tilde{\rho}_0\right\|_{L_x^q}\left\|\nabla \tilde{u}\right\|_{L_t^1 L_x^p}\right) \\
            & = 4\delta_0 \exp\left(\frac{C_{p,q}}{\kappa}\left\|\rho_0\right\|_{L_x^q}\left\|\nabla u\right\|_{L_t^1 L_x^p}\right)
        \end{align*}
        then fixing a mollifier $\tilde{\varphi}(x)\in\mathcal{S}(\mathbb{R}^d;\mathbb{R})$ satisfying the frequency cutoff \eqref{eq:frequencycutoff} as in the statement of Theorem \ref{thm:quantitative2}, one has
        \begin{align*}
            \left\|\rho_0 * \tilde{\varphi}_{\delta'}\right\|_{L_x^r} & \le\left\|\tilde{\rho} * \tilde{\varphi}_{\delta'}\right\|_{L_t^\infty L_x^r} \\
            & \le \kappa + A\left\|\tilde{\rho}_0*\tilde{\varphi}_{\delta_0}\right\|_{L_x^r} \\
            & \le \kappa + A\left\|\rho(\cdot,T)*\tilde{\varphi}_{\delta_0}\right\|_{L_x^r},
        \end{align*}
        with constants $C_{p,q} > 0$ depending only on the parameters $1<p,q\le\infty$ and the dimension $d\ge1$, and $A>0$ depending only on the dimension $d\ge1$.

        Meanwhile,
        \begin{align*}
            & \left\|\rho(\cdot,T)*\tilde{\varphi}_{\delta_0}\right\|_{L_x^r} \\
            & \qquad = \sup_{\left\|\phi\right\|_{L_x^{r'}}\le1}\int_{\mathbb{R}^d}\int_{\mathbb{R}^d}\phi(y)\tilde{\varphi}_{\delta_0}(y-x)\rho(x,T)\;dy\;dx \\
            & \qquad \le \left\|\rho(\cdot,T)\right\|_{W_x^{-1,r}}\sup_{\left\|\phi\right\|_{L_x^{r'}}\le1}\left(\left\|\int_{\mathbb{R}^d}\phi(y)\tilde{\varphi}_{\delta_0}(y-x)\;dy\right\|_{L_x^{r'}}+\left\|\nabla\int_{\mathbb{R}^d}\phi(y)\tilde{\varphi}_{\delta_0}(y-x)\;dy\right\|_{L_x^{r'}}\right) \\
            & \qquad \le \left\|\rho(\cdot,T)\right\|_{W_x^{-1,r}}\left(\left\|\tilde{\varphi}_{\delta_0}\right\|_{L_x^1}+\left\|\nabla\tilde{\varphi}_{\delta_0}\right\|_{L_x^1}\right) \\
            & \qquad \le A' \left(1+\frac{1}{\delta_0}\right)\left\|\rho(\cdot,T)\right\|_{W_x^{-1,r}},
        \end{align*}
        for a constant $A'>0$ depending only on the choice of mollifier $\tilde{\varphi}(x)\in\mathcal{S}(\mathbb{R}^d;\mathbb{R})$ and the dimension $d\ge1$.

        Combining the above, we see that
        \begin{align*}
            \left\|\rho(\cdot,T)\right\|_{W_x^{-1,r}} & \ge \frac{1}{AA'}\left(1+\frac{1}{\delta_0}\right)^{-1}\left(\left\|\rho_0*\tilde{\varphi}_{\delta'}\right\|_{L_x^r}-\kappa\right) \\
            & \ge \frac{1}{AA'}\left(1+\frac{4}{\delta'}\exp\left(\frac{C_{p,q}}{\kappa}\frac{\left\|\rho_0\right\|_{L_x^q}}{\left\|\rho_0\right\|_{L_x^r}}\left\|\nabla u\right\|_{L_t^1 L_x^p}\right)\right)^{-1}\left(\left\|\rho_0*\tilde{\varphi}_{\delta'}\right\|_{L_x^r}-\kappa\right),
        \end{align*}
        for any choice of $\delta'>0$. In particular, if $\delta'\le 1$ is such that
        \begin{equation}\label{eq:initialgeometricmixing}
            \left\|\rho_0*\tilde{\varphi}_{\delta'}\right\|_{L_x^r} \ge \frac{1}{4} \left\|\rho_0\right\|_{L_x^r},
        \end{equation}
        then taking $\kappa=\frac{1}{8}\left\|\rho_0\right\|_{L_x^r}$ we have the bound
        \begin{equation*}
            \left\|\rho(\cdot,T)\right\|_{W_x^{-1,r}} \ge \frac{\delta'\left\|\rho_0\right\|_{L_x^r}}{40AA'}\exp\left(-8C_{p,q}\frac{\left\|\rho_0\right\|_{L_x^q}}{\left\|\rho_0\right\|_{L_x^r}}\left\|\nabla u\right\|_{L_t^1 L_x^p}\right).
        \end{equation*}

        We now aim to generalise the condition \eqref{eq:initialgeometricmixing} to any positive mollifier $\varphi(x)\in L^1(\mathbb{R}^d;\mathbb{R})$, $\varphi(x)\ge0$, $\left\|\varphi\right\|_{L_x^1}=1$. Let $0<\delta\le 1$ be small enough that
        \begin{equation*}
            \left\|\rho_0*\varphi_\delta\right\|_{L_x^r}\ge\frac{1}{2}\left\|\rho_0\right\|_{L_x^r}.
        \end{equation*}

        Since $(\varphi * \tilde{\varphi}_\alpha)(x) \xrightarrow{\alpha\to0} \varphi(x)$ in $L_x^1$, there exists some $0<\alpha\le1$ such that $\left\|\varphi*\tilde{\varphi}_\alpha-\varphi\right\|_{L_x^1}\le\frac{1}{4}$, and so for all $\delta>0$, also
        \begin{equation*}
            \left\|\varphi_\delta*\tilde{\varphi}_{\alpha\delta}-\varphi_\delta\right\|_{L_x^1}\le\frac{1}{4}.
        \end{equation*}

        Let $\delta'=\alpha\delta\le 1$. Then since $\left\|\varphi\right\|_{L_x^1}=1$,
        \begin{align*}
            \left\|\rho_0*\tilde{\varphi}_{\delta'}\right\|_{L_x^r} & \ge \left\|\rho_0*\varphi_\delta*\tilde{\varphi}_{\delta'}\right\|_{L_x^r} \\
            & \ge \left\|\rho_0*\varphi_\delta\right\|_{L_x^r} - \left\|\rho_0*\varphi_\delta*\tilde{\varphi}_{\delta'}-\rho_0*\varphi_\delta\right\|_{L_x^r} \\
            & \ge \frac{1}{2}\left\|\rho_0\right\|_{L_x^r} - \left\|\rho_0\right\|_{L_x^r}\left\|\varphi_\delta*\tilde{\varphi}_{\alpha\delta}-\varphi_\delta\right\|_{L_x^1} \\
            & \ge \frac{1}{4}\left\|\rho_0\right\|_{L_x^r},
        \end{align*}
        and so
        \begin{equation*}
            \left\|\rho(\cdot,T)\right\|_{W_x^{-1,r}} \ge \frac{\alpha\delta\left\|\rho_0\right\|_{L_x^r}}{40AA'}\exp\left(-8C_{p,q}\frac{\left\|\rho_0\right\|_{L_x^q}}{\left\|\rho_0\right\|_{L_x^r}}\left\|\nabla u\right\|_{L_t^1 L_x^p}\right),
        \end{equation*}
        where $\alpha>0$ depends only on the choice of positive mollifier $\varphi(x)\in L^1(\mathbb{R}^d;\mathbb{R})$.
    \end{proof}

    \begin{remark}
        We remark that the dependence of the decay rate on the conserved ratio $\frac{\left\|\rho_0\right\|_{L_x^q}}{\left\|\rho_0\right\|_{L_x^r}}$ is necessary by dimensional analysis. This ratio captures the relevant length scale of the passive scalar and ensures that the full rate
        \begin{equation*}
            \frac{\left\|\rho_0\right\|_{L_x^q}}{\left\|\rho_0\right\|_{L_x^r}}\left\|\nabla u\right\|_{L_t^1L_x^p},
        \end{equation*}
        is dimensionless. In the case $p=\infty$ this reduces to the classical decay rate $\left\|\nabla u\right\|_{L_t^1L_x^\infty}$, however for $p<\infty$ it is important to observe that $\left\|\nabla u\right\|_{L_t^1L_x^p}$ has units of length to the power $\frac{d}{p}$ and so cannot bound the mixing rate by its self.
    \end{remark}

    \medskip

    We next apply Lemma \ref{thm:calderon} to show that the passive scalar propagates regularity corresponding to a logarithm of a derivative, also studied in \cite{ben2019convergence,leger2018new,brue2021sharp,brue2021advection,meyer2022propagation}:

    \begin{proposition}[Propagation of regularity]\label{thm:propagation}
        Consider a weak solution $\rho(x,t):\mathbb{R}^d\times[0,T]\to\mathbb{R}$, with $\rho(x,t)\in L_t^\infty L_x^q$, and initial datum $\rho_0(x)\in L_x^q$, to the transport equation
        \begin{equation*}
            \begin{cases}
                \frac{\partial \rho}{\partial t}(x,t)+\nabla \cdot(u(x,t)\rho(x,t))=0, \\
                \rho(x,0)=\rho_0(x),
            \end{cases}
        \end{equation*}
        along a divergence-free vector field $u(x,t):\mathbb{R}^d\times[0,T]\to\mathbb{R}^d$, with $\nabla u(x,t)\in L_t^1L_x^p$, for $\frac{1}{p}+\frac{1}{q}=\frac{1}{r}$, with $1<p,q\le\infty$, $1\le r<\infty$.
        
        Then for the convolution kernel $K(x)=\log(1-\Delta)$, given by its Fourier transform
        \begin{equation*}
            \hat{K}(\xi)=\log(1+\left|\xi\right|^2),
        \end{equation*}
        one has propagation of
        \begin{equation*}
            \frac{d}{dt}\left\|(\rho*K)(\cdot,t)\right\|_{L_x^r}\le C_{p,q} \left\|\rho(\cdot,t)\right\|_{L_x^q}\left\|\nabla u(\cdot,t)\right\|_{L_x^p},
        \end{equation*}
        where the constant $C_{p,q}>0$ depends only on the parameters $1<p,q\le\infty$ and the dimension $d\ge1$.
    \end{proposition}
    \begin{proof}
        We apply Lemma \ref{thm:calderon} with the kernel $K(x)=\log(1-\Delta)$, whose Fourier transform $\hat{K}(\xi):\mathbb{R}^d\to\mathbb{C}$ is given by $\hat{K}(\xi)=\log\left(1+|\xi|^2\right)$. It remains to show that the kernel
        \begin{equation*}
            K_{i,j}(x) = \frac{\partial(x_j K)}{\partial x_i}(x),
        \end{equation*}
        whose Fourier transform is given by
        \begin{align*}
            \hat{K}_{i,j}(\xi) & = -\xi_i\frac{\partial}{\partial \xi_j}\hat{K}(\xi) \\
            & = -\frac{2\xi_i\xi_j}{1+|\xi|^2},
        \end{align*}
        has finite Calderon-Zygmund norm for all $i,j=1,\dots,d$.

        To this end, we apply the following version of the Mikhlin-multiplier theorem: \cite{muscalu2013classical}*{Theorem 8.2}. We must show that
        \begin{equation*}
            |\xi|^n |\nabla^n\hat{K}_{i,j}(\xi)|,
        \end{equation*}
        is bounded for all $n=0,\dots,d+2$. Indeed, we have, by direct calculation, the bound
        \begin{align*}
            |\nabla^n\hat{K}_{i,j}(\xi)| & \le C_n \frac{1+|\xi|^{n+2}}{(1+|\xi|^2)^{n+1}},
        \end{align*}
        for all $n\ge0$, and some $C_n>0$.
    \end{proof}

    \medskip

    We next apply Theorem \ref{thm:quantitative2} to give a quantitative convergence rate for the vanishing diffusivity/viscosity limit of passive scalar transport, also studied in \cite{seis2018optimal,brue2021advection}:
    
    \begin{proposition}[Convergence of the vanishing diffusion limit]\label{thm:vanishingviscosity}
        Consider a weak solution $\rho(x,t):\mathbb{R}^d\times[0,T]\to\mathbb{R}$, with $\rho(x,t)\in L_t^\infty L_x^q$, and initial datum $\rho_0(x)\in L_x^q$, to the transport equation
        \begin{equation*}
            \begin{cases}
                \frac{\partial \rho}{\partial t}(x,t)+\nabla \cdot(u(x,t)\rho(x,t))=0, \\
                \rho(x,0)=\rho_0(x),
            \end{cases}
        \end{equation*}
        along a divergence-free vector field $u(x,t):\mathbb{R}^d\times[0,T]\to\mathbb{R}^d$, with $\nabla u(x,t)\in L_t^1L_x^p$, for $\frac{1}{p}+\frac{1}{q}=\frac{1}{r}$, with $1<p,q\le\infty$, $1\le r<\infty$.

        For any $\nu>0$ let $\rho_\nu(x,t)\in L_t^\infty L_x^q$ be a weak solution to the transport-diffusion equation
        \begin{equation*}
            \begin{cases}
                \frac{\partial \rho_\nu}{\partial t}(x,t)+\nabla \cdot(u(x,t)\rho_\nu(x,t))-\nu\Delta\rho_\nu(x,t)=0, \\
                \rho_\nu(x,0)=\rho_0(x),
            \end{cases}
        \end{equation*}
        with the same initial datum $\rho_0(x)\in L_x^q$.

        Fix a mollifier/frequency cutoff $\varphi(x)\in\mathcal{S}(\mathbb{R}^d;\mathbb{R})$ as in Definition \ref{def:frequencycutoff}.
        
        We have the following weak convergence of $\rho_\nu(x,t)\xrightharpoonup{}\rho(x,t)$ as $\nu\to0$:

        For any $\delta>0$, if $\nu T \le \delta^2 \frac{\left\|\rho_0\right\|_{L_x^q}}{\left\|\rho_0\right\|_{L_x^r}}\left\|\nabla u\right\|_{L_t^1 L_x^p}$, then
        \begin{equation*}
            \left\|(\rho-\rho_\nu)*\varphi_\delta\right\|_{L_t^\infty L_x^r} \le C_{p,q} \left\|\rho_0\right\|_{L_x^q}\left\|\nabla u\right\|_{L_t^1 L_x^p}\left(\log\left(\frac{\delta^2}{\nu T}\left(\frac{\left\|\rho_0\right\|_{L_x^q}}{\left\|\rho_0\right\|_{L_x^r}}\left\|\nabla u\right\|_{L_t^1 L_x^p}\right)\right)\right)^{-1},
        \end{equation*}
        where the constant $C_{p,q}>0$ depends only on the parameters $1<p,q\le\infty$, the choice of mollifier $\varphi(x)\in\mathcal{S}(\mathbb{R}^d;\mathbb{R})$, and the dimension $d\ge1$.
    \end{proposition}
    \begin{corollary}\label{cor:vanishingviscosity}
        If
        \begin{equation*}
            \nu T \le \min\left(1,\bigg(\frac{\left\|\rho_0\right\|_{L_x^q}}{\left\|\rho_0\right\|_{L_x^r}}\left\|\nabla u\right\|_{L_t^1 L_x^p}\bigg)^6\right),
        \end{equation*}
        then
        \begin{equation*}
            \left\|\rho-\rho_\nu\right\|_{L_t^\infty W_x^{-1,r}} \le C_{p,q} \left\|\rho_0\right\|_{L_x^q}\left\|\nabla u\right\|_{L_t^1 L_x^p}\left(\log\left(\frac{1}{\nu T}\right)\right)^{-1},
        \end{equation*}
        where the constant $C_{p,q}>0$ depends only on the parameters $1<p,q\le\infty$ and the dimension $d\ge1$.
    \end{corollary}
    \begin{proof}[Proof of Proposition \ref{thm:vanishingviscosity}]
        Fix a mollifier $\varphi(x)\in\mathcal{S}(\mathbb{R}^d;\mathbb{R})$ satisfying the frequency-cutoff \eqref{eq:frequencycutoff}.
        
        We begin by writing $(\rho-\rho_\nu)(x,t)\in L_t^\infty L_x^q$ as a weak solution with zero initial datum to the following distributional forced transport equation along $u(x,t)$:
        \begin{equation*}
            \frac{\partial (\rho-\rho_\nu)}{\partial t}(x,t) + \nabla\cdot(u(x,t)(\rho-\rho_\nu)(x,t))=f(x,t),
        \end{equation*}
        with the distributional force given by
        \begin{equation*}
            f(x,t) = -\nu\Delta\rho_\nu(x,t).
        \end{equation*}

        By Theorem \ref{thm:quantitative2}, for any $\kappa > 0$, we have the following quantitative bound:

        For
        \begin{equation*}
            \delta' = \frac{\delta}{4}\exp\left(-\frac{C_{p,q}}{\kappa}\left\|\rho-\rho_\nu\right\|_{L_t^\infty L_x^q}\left\|\nabla u\right\|_{L_t^1 L_x^p}\right),
        \end{equation*}
        then
        \begin{align}
            & \left\|(\rho-\rho_\nu)*\varphi_\delta\right\|_{L_t^\infty L_x^r} \nonumber \\
            & \qquad \le \kappa + A\left\|f*\varphi_{\delta'}\right\|_{L_t^1L_x^r} \nonumber \\
            & \qquad \le \kappa + \nu A\left\|\Delta\varphi_{\delta'}\right\|_{L_x^1}\left\|\rho_\nu\right\|_{L_t^1L_x^r} \nonumber \\
            & \qquad \le \kappa + \nu A\delta'^{-2}\left\|\nabla\varphi\right\|_{L_x^1}T\left\|\rho_0\right\|_{L_x^{r}} \nonumber \\
            & \qquad \le \kappa + \nu \frac{16A}{\delta^{2}}\exp\left(2\frac{C_{p,q}}{\kappa}\left\|\rho-\rho_\nu\right\|_{L_t^\infty L_x^q}\left\|\nabla u\right\|_{L_t^1 L_x^p}\right)\left\|\nabla\varphi\right\|_{L_x^1}T\left\|\rho_0\right\|_{L_x^{r}} \nonumber \\
            & \qquad \le \kappa + \nu A'\delta^{-2}T\left\|\rho_0\right\|_{L_x^{r}}\exp\left(\frac{C_{p,q}'}{\kappa}\left\|\rho-\rho_\nu\right\|_{L_t^\infty L_x^q}\left\|\nabla u\right\|_{L_t^1 L_x^p}\right) \nonumber \\
            & \qquad \le \kappa + \nu A'\delta^{-2}T\left\|\rho_0\right\|_{L_x^{r}}\exp\left(\frac{C_{p,q}'}{\kappa}\left(\left\|\rho\right\|_{L_t^\infty L_x^q}+\left\|\rho_\nu\right\|_{L_t^\infty L_x^q}\right)\left\|\nabla u\right\|_{L_t^1 L_x^p}\right) \nonumber \\ 
            & \qquad \le \kappa + \nu A'\delta^{-2}T\left\|\rho_0\right\|_{L_x^{r}}\exp\left(\frac{C_{p,q}'}{\kappa}2\left\|\rho_0\right\|_{L_x^q}\left\|\nabla u\right\|_{L_t^1 L_x^p}\right), \label{eq:logarithmdecay1}
        \end{align}
        where the constants $C_{p,q}, C_{p,q}' > 0$ depends only on the parameters $1<p,q\le\infty$, the choice of mollifier $\varphi(x)\in\mathcal{S}(\mathbb{R}^d;\mathbb{R})$, and the dimension $d\ge1$, while $A,A'>0$ depend only on the choice of mollifier $\varphi(x)\in\mathcal{S}(\mathbb{R}^d;\mathbb{R})$ and the dimension $d\ge1$.

        Let
        \begin{gather*}
            \gamma = 2\left\|\rho_0\right\|_{L_x^q}\left\|\nabla u\right\|_{L_t^1 L_x^p},
        \end{gather*}
        and rewrite \eqref{eq:logarithmdecay1} as, for any $\kappa>0$,
        \begin{equation*}
            \left\|(\rho-\rho_\nu)*\varphi_\delta\right\|_{L_t^\infty L_x^r} \le \kappa + \nu T A'\left(2\delta^{2}\frac{\left\|\rho_0\right\|_{L_x^q}}{\left\|\rho_0\right\|_{L_x^r}}\left\|\nabla u\right\|_{L_t^1 L_x^p}\right)^{-1} \gamma\exp\left(C_{p,q}'{\frac{\gamma}{\kappa}}\right).
        \end{equation*}

        We now crudely optimise over $\kappa>0$:
        
        For $\nu T\le\delta^{2}\frac{\left\|\rho_0\right\|_{L_x^q}}{\left\|\rho_0\right\|_{L_x^r}}\left\|\nabla u\right\|_{L_t^1 L_x^p}$, let
        \begin{equation*}
            \kappa=2C_{p,q}'\gamma\left(\log\left(\frac{\delta^{2}}{\nu T}\frac{\left\|\rho_0\right\|_{L_x^q}}{\left\|\rho_0\right\|_{L_x^r}}\left\|\nabla u\right\|_{L_t^1 L_x^p}\right)\right)^{-1},
        \end{equation*}
        then, using $x^\frac{1}{2}\le \frac{2e^{-1}}{-\log x}$ for $x\in(0,1]$,
        \begin{align*}
            & \kappa + \nu TA'\left(2\delta^{2}\frac{\left\|\rho_0\right\|_{L_x^q}}{\left\|\rho_0\right\|_{L_x^r}}\left\|\nabla u\right\|_{L_t^1 L_x^p}\right)^{-1}\gamma \exp\left(C_{p,q}'\frac{\gamma}{\kappa}\right) \\
            & \qquad = 2C_{p,q}'\gamma\left(\log\left(\frac{\delta^{2}}{\nu T}\frac{\left\|\rho_0\right\|_{L_x^q}}{\left\|\rho_0\right\|_{L_x^r}}\left\|\nabla u\right\|_{L_t^1 L_x^p}\right)\right)^{-1} + \frac{1}{2}A'\gamma\left(\frac{\delta^2}{\nu T}\frac{\left\|\rho_0\right\|_{L_x^q}}{\left\|\rho_0\right\|_{L_x^r}}\left\|\nabla u\right\|_{L_t^1 L_x^p}\right)^{-\frac{1}{2}} \\
            & \qquad \le \left(2C_{p,q}'\gamma + A'e^{-1}\gamma\right)\left(\log\left(\frac{\delta^{2}}{\nu T}\frac{\left\|\rho_0\right\|_{L_x^q}}{\left\|\rho_0\right\|_{L_x^r}}\left\|\nabla u\right\|_{L_t^1 L_x^p}\right)\right)^{-1}.
        \end{align*}

        This bound concludes the proof of Proposition \ref{thm:vanishingviscosity}.
    \end{proof}
    \begin{proof}[Proof of Corollary \ref{cor:vanishingviscosity}]
        Continuing on from the proof of Proposition \ref{thm:vanishingviscosity}.

        In general, we may bound the weak norm $W_x^{-1,r}$ by
        \begin{align*}
            \left\|f\right\|_{W_x^{-1,r}} & = \sup_{\left\|\phi\right\|_{L_x^{r'}}+\left\|\nabla\phi\right\|_{L_x^{r'}}\le1}\int_{\mathbb{R}^d}f(x)\phi(x)\;dx \\
            & = \sup_{\left\|\phi\right\|_{L_x^{r'}}+\left\|\nabla\phi\right\|_{L_x^{r'}}\le1}\left(\int_{\mathbb{R}^d}(f*\varphi_\delta)(x)\phi(x)\;dx + \int_{\mathbb{R}^d}(f-f*\varphi_\delta)(x)\phi(x)\;dx\right) \\
            & \le \sup_{\left\|\phi\right\|_{L_x^{r'}}+\left\|\nabla\phi\right\|_{L_x^{r'}}\le1}\left(\left\|f*\varphi_\delta\right\|_{L_x^r}\left\|\phi\right\|_{L_x^{r'}} + \delta \left\|f\right\|_{L_x^r}\left\|\nabla\phi\right\|_{L_x^{r'}}\right) \\ \\
            & \lesssim \left\|f*\varphi_\delta\right\|_{L_x^r} + \delta \left\|f\right\|_{L_x^r}.
        \end{align*}


        Applying this to $(\rho-\rho_\nu)(x,t)$, and the previous bounds, gives
        \begin{align}
            \left\|\rho-\rho_\nu\right\|_{L_t^\infty W_x^{-1,r}} & \lesssim \left\|\rho_0\right\|_{L_x^q}\left\|\nabla u\right\|_{L_t^1 L_x^p}\left(\log\left(\frac{\delta^2}{\nu T}\frac{\left\|\rho_0\right\|_{L_x^q}}{\left\|\rho_0\right\|_{L_x^r}}\left\|\nabla u\right\|_{L_t^1 L_x^p}\right)\right)^{-1} + \delta\left\|\rho-\rho_\nu\right\|_{L_x^r} \nonumber \\
            & \lesssim \left\|\rho_0\right\|_{L_x^q}\left\|\nabla u\right\|_{L_t^1 L_x^p}\left(\log\left(\frac{\delta^2}{\nu T}\frac{\left\|\rho_0\right\|_{L_x^q}}{\left\|\rho_0\right\|_{L_x^r}}\left\|\nabla u\right\|_{L_t^1 L_x^p}\right)\right)^{-1} + \delta\left\|\rho_0\right\|_{L_x^r}. \label{eq:optimisedelta2}
        \end{align}

        We now crudely optimise over $\delta>0$:

        For $\nu T\le \left(\frac{\left\|\rho_0\right\|_{L_x^q}}{\left\|\rho_0\right\|_{L_x^r}}\left\|\nabla u\right\|_{L_t^1 L_x^p}\right)^3$ let
        \begin{equation*}
            \delta = \frac{\left\|\rho_0\right\|_{L_x^q}}{\left\|\rho_0\right\|_{L_x^r}}\left\|\nabla u\right\|_{L_t^1 L_x^p}\left(\log\left(\frac{1}{\nu T}\left(\frac{\left\|\rho_0\right\|_{L_x^q}}{\left\|\rho_0\right\|_{L_x^r}}\left\|\nabla u\right\|_{L_t^1 L_x^p}\right)^3\right)\right)^{-1},
        \end{equation*}
        and so, using that $-\log\log x \ge -e^{-1}\log x$ for $x\ge1$,
        \begin{align*}
            \log\left(\frac{\delta^2}{\nu T}\frac{\left\|\rho_0\right\|_{L_x^q}}{\left\|\rho_0\right\|_{L_x^r}}\left\|\nabla u\right\|_{L_t^1 L_x^p}\right) & = \log\left(\frac{1}{\nu T}\left(\frac{\left\|\rho_0\right\|_{L_x^q}}{\left\|\rho_0\right\|_{L_x^r}}\left\|\nabla u\right\|_{L_t^1 L_x^p}\right)^3\right) - 2 \log\log\left(\frac{1}{\nu T}\left(\frac{\left\|\rho_0\right\|_{L_x^q}}{\left\|\rho_0\right\|_{L_x^r}}\left\|\nabla u\right\|_{L_t^1 L_x^p}\right)^3\right) \\
            & \ge (1-2e^{-1})\log\left(\frac{1}{\nu T}\left(\frac{\left\|\rho_0\right\|_{L_x^q}}{\left\|\rho_0\right\|_{L_x^r}}\left\|\nabla u\right\|_{L_t^1 L_x^p}\right)^3\right),
        \end{align*}

        Therefore \eqref{eq:optimisedelta2} becomes
        \begin{equation*}
            \left\|\rho-\rho_\nu\right\|_{L_t^\infty W_x^{-1,r}} \lesssim \left\|\rho_0\right\|_{L_x^q}\left\|\nabla u\right\|_{L_t^1 L_x^p}\left(\log\left(\frac{1}{\nu T}\right)+\log\left(\frac{\left\|\rho_0\right\|_{L_x^q}}{\left\|\rho_0\right\|_{L_x^r}}\left\|\nabla u\right\|_{L_t^1 L_x^p}\right)^3\right)^{-1}.
        \end{equation*}

        Using now that for positive reals $a,b>0$, if $a\le\min(1,b^2)$, then $\left(\log\frac{b}{a}\right)^{-1}\le 2\left(\log\frac{1}{a}\right)^{-1}$, then if in addition
        \begin{equation*}
            \nu T \le \min\left(1,\left(\frac{\left\|\rho_0\right\|_{L_x^q}}{\left\|\rho_0\right\|_{L_x^r}}\left\|\nabla u\right\|_{L_t^1 L_x^p}\right)^6\right),
        \end{equation*}
        then
        \begin{equation*}
            \left(\log\left(\frac{1}{\nu T}\right)+\log\left(\frac{\left\|\rho_0\right\|_{L_x^q}}{\left\|\rho_0\right\|_{L_x^r}}\left\|\nabla u\right\|_{L_t^1 L_x^p}\right)^3\right)^{-1} \le 2\left(\log\left(\frac{1}{\nu T}\right)\right)^{-1},
        \end{equation*}
        which gives the result.
    \end{proof}

    \medskip

    Finally, we apply Theorem \ref{thm:quantitative2} to give a quantitative convergence rate for the passive scalar under any approximation of the vector field, also studied in \cite{seis2017quantitative,seis2018optimal}:
    \begin{proposition}[Approximation of the vector field]\label{thm:stability}
        Consider a weak solution $\rho(x,t):\mathbb{R}^d\times[0,T]\to\mathbb{R}$, with $\rho(x,t)\in L_t^\infty L_x^q$, and initial datum $\rho_0(x)\in L_x^q$, to the transport equation
        \begin{equation*}
            \begin{cases}
                \frac{\partial \rho}{\partial t}(x,t)+\nabla \cdot(u(x,t)\rho(x,t))=0, \\
                \rho(x,0)=\rho_0(x),
            \end{cases}
        \end{equation*}
        along a divergence-free vector field $u(x,t):\mathbb{R}^d\times[0,T]\to\mathbb{R}^d$, with $\nabla u(x,t)\in L_t^1L_x^p$, for $\frac{1}{p}+\frac{1}{q}=\frac{1}{r}$, with $1<p,q\le\infty$, $1\le r<\infty$.

        Let $\bar{\rho}(x,t)\in L_t^\infty L_x^q$ be a weak solution with the same initial datum $\rho_0(x)\in L_x^q$ along a measurable divergence-free vector field $\bar{u} \in L_t^1 L_x^p$:
        \begin{equation*}
            \begin{cases}
                \frac{\partial \bar{\rho}}{\partial t}(x,t)+\nabla \cdot(\bar{u}(x,t)\bar{\rho}(x,t))=0, \\
                \bar{\rho}(x,0)=\rho_0(x),
            \end{cases}
        \end{equation*}

        Fix a mollifier/frequency cutoff $\varphi(x)\in\mathcal{S}(\mathbb{R}^d;\mathbb{R})$ as in Definition \ref{def:frequencycutoff}.
        
        We have the following weak convergence of $\bar{\rho}(x,t)\xrightharpoonup{}\rho(x,t)$ as $\left\|u-\bar{u}\right\|_{L_t^1L_x^p}\to0$:
        
        For any $\delta>0$, if $\left\|u-\bar{u}\right\|_{L_t^1L_x^p}\le\delta\left\|\nabla u\right\|_{L_t^1 L_x^p}$, then
        \begin{equation*}
            \left\|(\rho-\bar{\rho})*\varphi_\delta\right\|_{L_t^\infty L_x^r} \le C_{p,q}\left\|\bar{\rho}\right\|_{L_t^\infty L_x^q}\left\|\nabla u\right\|_{L_t^1L_x^p}\left(\log\left(\frac{\delta\left\|\nabla u\right\|_{L_t^1 L_x^p}}{\left\|u-\bar{u}\right\|_{L_t^1L_x^p}}\right)\right)^{-1},
        \end{equation*}
        where the constant $C_{p,q}>0$ depends only on the parameters $1<p,q\le\infty$, the choice of mollifier $\varphi(x)\in\mathcal{S}(\mathbb{R}^d;\mathbb{R})$, and the dimension $d\ge1$.
    \end{proposition}
    \begin{corollary}\label{cor:stability}
        If
        \begin{equation*}
            \left\|u-\bar{u}\right\|_{L_t^1L_x^p}\le \left\|\nabla u\right\|_{L_t^1 L_x^p}\min\bigg(1,\bigg(\frac{\left\|\bar{\rho}\right\|_{L_t^\infty L_x^q}}{\left\|\bar{\rho}\right\|_{L_t^\infty L_x^r}}\left\|\nabla u\right\|_{L_t^1 L_x^p}\bigg)^2\bigg),
        \end{equation*}
        then
        \begin{equation*}
            \left\|\rho-\bar{\rho}\right\|_{L_t^\infty W_x^{-1,r}} \le C_{p,q} \left\|\bar{\rho}\right\|_{L_t^\infty L_x^q}\left\|\nabla u\right\|_{L_t^1 L_x^p}\left(\log\left(\frac{\left\|\nabla u\right\|_{L_t^1 L_x^p}}{\left\|u-\bar{u}\right\|_{L_t^1L_x^p}}\right)\right)^{-1},
        \end{equation*}
        where the constant $C_{p,q}>0$ depends only on the parameters $1<p,q\le\infty$ and the dimension $d\ge1$.
    \end{corollary}
    \begin{proof}[Proof of Proposition \ref{thm:stability}]
        We begin by writing $(\rho-\bar{\rho})(x,t)\in L_t^\infty L_x^q$ as a weak solution with zero initial datum to the following forced transport equation along $u(x,t)$:
        \begin{equation*}
            \frac{\partial (\rho-\bar{\rho})}{\partial t}(x,t) + \nabla\cdot(u(x,t)(\rho-\bar{\rho})(x,t))=f(x,t),
        \end{equation*}
        with the distributional force given by
        \begin{equation*}
            f(x,t) = -\nabla\cdot((u\bar{\rho}-\bar{u}\bar{\rho})(x,t)).
        \end{equation*}

        By Theorem \ref{thm:quantitative2}, for any $\kappa > 0$, we have the following quantitative bound:

        For
        \begin{equation*}
            \delta' = \frac{\delta}{4}\exp\left(-\frac{C_{p,q}}{\kappa}\left\|\rho-\bar{\rho}\right\|_{L_t^\infty L_x^q}\left\|\nabla u\right\|_{L_t^1 L_x^p}\right),
        \end{equation*}
        then
        \begin{align}
            & \left\|(\rho-\bar{\rho})*\varphi_\delta\right\|_{L_t^\infty L_x^r} \nonumber \\
            & \qquad \le \kappa + A\left\|f*\varphi_{\delta'}\right\|_{L_t^1L_x^r} \nonumber \\
            & \qquad \le \kappa + A\left\|\nabla\varphi_{\delta'}\right\|_{L_x^1}\left\|u\bar{\rho}-\bar{u}\bar{\rho}\right\|_{L_t^1L_x^r} \nonumber \\
            & \qquad \le \kappa + A\delta'^{-1}\left\|\nabla\varphi\right\|_{L_x^1}\left\|u-\bar{u}\right\|_{L_t^1L_x^p}\left\|\bar{\rho}\right\|_{L_t^\infty L_x^q} \nonumber \\
            & \qquad \le \kappa + \frac{4A}{\delta}\exp\left(\frac{C_{p,q}}{\kappa}\left\|\rho-\bar{\rho}\right\|_{L_t^\infty L_x^q}\left\|\nabla u\right\|_{L_t^1 L_x^p}\right)\left\|u-\bar{u}\right\|_{L_t^1L_x^p}\left\|\bar{\rho}\right\|_{L_t^\infty L_x^q} \nonumber \\
            & \qquad \le \kappa + A'\delta^{-1}\left\|u-\bar{u}\right\|_{L_t^1L_x^p}\left\|\bar{\rho}\right\|_{L_t^\infty L_x^q}\exp\left(\frac{C_{p,q}'}{\kappa}\left\|\bar{\rho}\right\|_{L_t^\infty L_x^q}\left\|\nabla u\right\|_{L_t^1 L_x^p}\right), \label{eq:logarithmdecay0}
        \end{align}
        where the constants $C_{p,q}, C_{p,q}' > 0$ depends only on the parameters $1<p,q\le\infty$, the choice of mollifier $\varphi(x)\in\mathcal{S}(\mathbb{R}^d;\mathbb{R})$, and the dimension $d\ge1$, while $A,A'>0$ depend only on the choice of mollifier $\varphi(x)\in\mathcal{S}(\mathbb{R}^d;\mathbb{R})$ and the dimension $d\ge1$.

        Let
        \begin{gather*}
            \gamma = 2\left\|\bar{\rho}\right\|_{L_t^\infty L_x^q}\left\|\nabla u\right\|_{L_t^1 L_x^p},
        \end{gather*}
        and rewrite \eqref{eq:logarithmdecay0} as, for any $\kappa>0$,
        \begin{equation*}
            \left\|(\rho-\bar{\rho})*\varphi_\delta\right\|_{L_t^\infty L_x^r} \le \kappa + A'\frac{\left\|u-\bar{u}\right\|_{L_t^1 L_x^p}}{2\delta\left\|\nabla u\right\|_{L_t^1 L_x^p}}\gamma \exp\left(C_{p,q}'{\frac{\gamma}{\kappa}}\right).
        \end{equation*}

        We now crudely optimise over $\kappa>0$:
        
        For $\left\|u-\bar{u}\right\|_{L_t^1 L_x^p}\le\delta\left\|\nabla u\right\|_{L_t^1 L_x^p}$, let
        \begin{equation*}
            \kappa=2C_{p,q}'\gamma\left(\log\frac{\delta\left\|\nabla u\right\|_{L_t^1 L_x^p}}{\left\|u-\bar{u}\right\|_{L_t^1 L_x^p}}\right)^{-1},
        \end{equation*}
        then, using $x^\frac{1}{2}\le \frac{2e^{-1}}{-\log x}$ for $x\in(0,1]$,
        \begin{align*}
            & \kappa + A'\frac{\left\|u-\bar{u}\right\|_{L_t^1 L_x^p}}{2\delta\left\|\nabla u\right\|_{L_t^1 L_x^p}}\gamma \exp\left(C_{p,q}'\frac{\gamma}{\kappa}\right) \\
            & \qquad = 2C_{p,q}'\gamma\left(\log\frac{\delta\left\|\nabla u\right\|_{L_t^1 L_x^p}}{\left\|u-\bar{u}\right\|_{L_t^1 L_x^p}}\right)^{-1} + \frac{1}{2}A'\gamma\left(\frac{\delta\left\|\nabla u\right\|_{L_t^1 L_x^p}}{\left\|u-\bar{u}\right\|_{L_t^1 L_x^p}}\right)^{-\frac{1}{2}} \\
            & \qquad \le \left(2C_{p,q}'\gamma+A'e^{-1}\gamma\right)\left(\log\frac{\delta\left\|\nabla u\right\|_{L_t^1 L_x^p}}{\left\|u-\bar{u}\right\|_{L_t^1 L_x^p}}\right)^{-1}.
        \end{align*}

        This bound concludes the proof of Proposition \ref{thm:stability}.
    \end{proof}
    \begin{proof}[Proof of Corollary \ref{cor:stability}]
        Continuing on from the proof of Proposition \ref{thm:stability}.


        As in the proof of Corollary \ref{cor:vanishingviscosity}, we may bound the weak norm by
        \begin{equation*}
            \left\|f\right\|_{W_x^{-1,r}} \lesssim \left\|f*\varphi_\delta\right\|_{L_x^r} + \delta \left\|f\right\|_{L_x^r}.
        \end{equation*}

        Applying this to $(\rho-\bar{\rho})(x,t)$, and the previous bounds, gives
        \begin{align}
            \left\|\rho-\bar{\rho}\right\|_{L_t^\infty W_x^{-1,r}} & \lesssim \left\|\bar{\rho}\right\|_{L_t^\infty L_x^q}\left\|\nabla u\right\|_{L_t^1 L_x^p}\left(\log\left(\frac{\delta\left\|\nabla u\right\|_{L_t^1 L_x^p}}{\left\|u-\bar{u}\right\|_{L_t^1L_x^p}}\right)\right)^{-1} + \delta\left\|\rho-\bar{\rho}\right\|_{L_x^r} \nonumber \\
            & \lesssim \left\|\bar{\rho}\right\|_{L_t^\infty L_x^q}\left\|\nabla u\right\|_{L_t^1 L_x^p}\left(\log\left(\frac{\delta\left\|\nabla u\right\|_{L_t^1 L_x^p}}{\left\|u-\bar{u}\right\|_{L_t^1L_x^p}}\right)\right)^{-1} + \delta\left\|\bar{\rho}\right\|_{L_t^\infty L_x^r}. \label{eq:optimisedelta1}
        \end{align}

        We now crudely optimise over $\delta>0$:

        For $\left\|u-\bar{u}\right\|_{L_t^1L_x^p}\le \delta\left\|\nabla u\right\|_{L_t^1 L_x^p}$ let
        \begin{equation*}
            \delta = \frac{\left\|\bar{\rho}\right\|_{L_t^\infty L_x^q}\left\|\nabla u\right\|_{L_t^1 L_x^p}}{\left\|\bar{\rho}\right\|_{L_t^\infty L_x^r}}\left(\log\left(\frac{\left\|\bar{\rho}\right\|_{L_t^\infty L_x^q}\left\|\nabla u\right\|_{L_t^1 L_x^p}^2}{\left\|u-\bar{u}\right\|_{L_t^1L_x^p}\left\|\bar{\rho}\right\|_{L_t^\infty L_x^r}}\right)\right)^{-1},
        \end{equation*}
        and so, using that $-\log\log x \ge -e^{-1}\log x$ for $x\ge1$,
        \begin{align*}
            \log\left(\frac{\delta\left\|\nabla u\right\|_{L_t^1 L_x^p}}{\left\|u-\bar{u}\right\|_{L_t^1L_x^p}}\right) & = \log\left(\frac{\left\|\bar{\rho}\right\|_{L_t^\infty L_x^q}\left\|\nabla u\right\|_{L_t^1 L_x^p}^2}{\left\|u-\bar{u}\right\|_{L_t^1L_x^p}\left\|\bar{\rho}\right\|_{L_t^\infty L_x^r}}\right) -  \log\log\left(\frac{\left\|\bar{\rho}\right\|_{L_t^\infty L_x^q}\left\|\nabla u\right\|_{L_t^1 L_x^p}^2}{\left\|u-\bar{u}\right\|_{L_t^1L_x^p}\left\|\bar{\rho}\right\|_{L_t^\infty L_x^r}}\right) \\
            & \ge (1-e^{-1})\log\left(\frac{\left\|\bar{\rho}\right\|_{L_t^\infty L_x^q}\left\|\nabla u\right\|_{L_t^1 L_x^p}^2}{\left\|u-\bar{u}\right\|_{L_t^1L_x^p}\left\|\bar{\rho}\right\|_{L_t^\infty L_x^r}}\right),
        \end{align*}

        Therefore, \eqref{eq:optimisedelta1} becomes
        \begin{equation*}
            \left\|\rho-\bar{\rho}\right\|_{L_t^\infty W_x^{-1,r}} \lesssim \left\|\bar{\rho}\right\|_{L_t^\infty L_x^q}\left\|\nabla u\right\|_{L_t^1 L_x^p}\left(\log\left(\frac{\left\|\nabla u\right\|_{L_t^1 L_x^p}}{\left\|u-\bar{u}\right\|_{L_t^1L_x^p}}\right)+ \log\left(\frac{\left\|\bar{\rho}\right\|_{L_t^\infty L_x^q}}{\left\|\bar{\rho}\right\|_{L_t^\infty L_x^r}}\left\|\nabla u\right\|_{L_t^1 L_x^p}\right)\right)^{-1}.
        \end{equation*}

        Using now that for positive reals $a,b>0$, if $a\le\min(1,b^2)$, then $\left(\log\frac{b}{a}\right)^{-1}\le 2\left(\log\frac{1}{a}\right)^{-1}$, we see that if in addition
        \begin{equation*}
            \frac{\left\|u-\bar{u}\right\|_{L_t^1L_x^p}}{\left\|\nabla u\right\|_{L_t^1 L_x^p}} \le \min\left(1,\left(\frac{\left\|\bar{\rho}\right\|_{L_t^\infty L_x^q}}{\left\|\bar{\rho}\right\|_{L_t^\infty L_x^r}}\left\|\nabla u\right\|_{L_t^1 L_x^p}\right)^2\right),
        \end{equation*}
        then
        \begin{equation*}
            \left(\log\left(\frac{\left\|\nabla u\right\|_{L_t^1 L_x^p}}{\left\|u-\bar{u}\right\|_{L_t^1L_x^p}}\right)+ \log\left(\frac{\left\|\bar{\rho}\right\|_{L_t^\infty L_x^q}}{\left\|\bar{\rho}\right\|_{L_t^\infty L_x^r}}\left\|\nabla u\right\|_{L_t^1 L_x^p}\right)\right)^{-1} \le 2\left(\log\left(\frac{\left\|\nabla u\right\|_{L_t^1 L_x^p}}{\left\|u-\bar{u}\right\|_{L_t^1L_x^p}}\right)\right)^{-1},
        \end{equation*}
        which gives the result.
    \end{proof}

    \begin{remark}
        We remark that the results of this section do not extend beyond the DiPerna-Lions regime $\frac{1}{p}+\frac{1}{q}\le1$. Outside this class of weak solutions $\rho(x,t)\in L_t^\infty L_x^q$ there are counterexamples to uniqueness (and thus any quantitative estimates), see \cite{modena2018non,modena2020convex,brue2024sharp}. We remark, however, that there is still a unique regular Lagrangian flow in this setting \cite{ambrosio2014continuity}, and whether solutions selected by this flow satisfy quantitative estimates analogous to these results is an open question.
    \end{remark}

\section{Estimates on the DiPerna-Lions commutator}\label{section4}

    Finally, we discuss how some of the results in this paper are related to decay rates for the standard DiPerna-Lions commutator in \cite{diperna1989ordinary}. 
    
    The DiPerna-Lions commutator controls the transport error of the mollification $(\rho*\varphi_\delta)(x)$:
    \begin{equation*}
        \frac{\partial}{\partial t}(\rho*\varphi_\delta)(x,t) + u(x,t)\cdot\nabla(\rho*\varphi_\delta)(x,t) = R_\delta(u,\rho)(x,t),
    \end{equation*}
    where we have kept the dependence of the commutator $R_\delta(u,\rho)(x,t)$ on the passive scalar $\rho(x,t)$ and the vector field $u(x,t)$. The DiPerna-Lions commutator is a special case of the commutator $R(x,t)$ in Lemma \ref{thm:calderon} by taking the kernel $K(x)=\varphi_\delta(x)$. Applying Lemma \ref{thm:calderon} as in the proof of Theorem \ref{thm:quantitative2}, the harmonic estimate in the proof of Theorem \ref{thm:quantitative2} can be rephrased as the following integral estimate on the decay of the DiPerna-Lions commutator. Suppressing the time-dependence for readability:
    \begin{proposition}[Integral estimate on the DiPerna-Lions commutator]\label{thm:expdecay}
        Consider a passive scalar $\rho(x):\mathbb{R}^d\to\mathbb{R}$, with $\rho(x)\in L_x^q$, and a divergence-free vector field $u(x):\mathbb{R}^d\to\mathbb{R}^d$, with $\nabla u(x)\in L_x^p$, for $\frac{1}{p}+\frac{1}{q}=\frac{1}{r}$, with $1<p,q\le\infty$, $1\le r<\infty$.

        Let $\varphi(x)\in \mathcal{S}(\mathbb{R}^d;\mathbb{R})$ (a Schwartz function) be a standard mollifier and denote by
        \begin{equation*}
            R_\delta(u,\rho)(x) = u(x)\cdot\nabla(\rho*\varphi_\delta)(x) - \nabla\cdot((u\rho)*\varphi_\delta)(x),
        \end{equation*}
        the DiPerna-Lions commutator, as in \cite{diperna1989ordinary}.
        
        Then for any $0<\delta_1<{\delta_2}$,
        \begin{equation*}
            \left\| \int_{\delta_1}^{\delta_2} R_\delta(u,\rho)(x) \; \frac{d\delta}{\delta} \right\|_{L_x^r} \le C_{p,q} \left\|\nabla u\right\|_{L_x^p}\left\|\rho\right\|_{L_x^q},
        \end{equation*}
        where the constant $C_{p,q} > 0$ depends only on the parameters $1<p,q\le\infty$, the choice of mollifier $\varphi(x)\in\mathcal{S}(\mathbb{R}^d;\mathbb{R})$, and the dimension $d\ge1$.
    \end{proposition}
    \begin{proof}
        Throughout, $\lesssim$ will denote less than or equal to up to a constant depending only on $p,q$, the dimension $d$, and the mollifier $\varphi(x)\in\mathcal{S}(\mathbb{R}^d;\mathbb{R})$ (and in particular not on $\delta_1,{\delta_2}$).

        We first rewrite
        \begin{align*}
            \int_{\delta_1}^{\delta_2} R_\delta(u,\rho)(x)\;\frac{d\delta}{\delta} & = \int_{\delta_1}^{\delta_2}u(x)\cdot(\rho*\nabla \varphi_\delta)(x) - \nabla\cdot((u\rho)*\varphi_\delta)(x)\;\frac{d\delta}{\delta} \\
            & = u(x)\cdot(\rho * K)(x) - \nabla\cdot((u\rho)*K)(x),
        \end{align*}
        in terms of the kernel
        \begin{equation*}
            K(x) = \int_{\delta_1}^{\delta_2}\varphi_\delta(x)\frac{d\delta}{\delta}.
        \end{equation*}
        
        We now apply Lemma \ref{thm:calderon} as in the proof of Theorem \ref{thm:quantitative2} to obtain the estimate
        \begin{equation}\label{eq:integralcommutator}
            \left\|\int_{\delta_1}^{\delta_2} R_\delta(u,\rho)(x)\;\frac{d\delta}{\delta}\right\|_{L_x^r} \le C_{p,q} \left\|\rho\right\|_{L_x^q}\left\|\nabla u\right\|_{L_x^p}\left(\sum_{i,j=1}^d\left\|\frac{\partial(x_j K)}{\partial x_i}\right\|_{CZ}\right),
        \end{equation}
        where the constant $C_{p,q}>0$ depends only on the parameters $1<p,q\le\infty$ and the dimension $d\ge1$, and where $\left\|\cdot\right\|_{CZ}$ refers to the Calderon-Zygmund norm, as defined in Lemma \ref{thm:calderon}.
        
        Now let
        \begin{align*}
            K_{i,j}'(x) & = \frac{\partial (x_j K)}{\partial x_i}(x) = \int_{\delta_1}^{\delta_2}\frac{\partial (x_j\varphi)}{\partial x_i}\left(\frac{x}{\delta}\right)\;\frac{d\delta}{\delta^{d+1}}.
        \end{align*}
        
        As in the proof of Theorem \ref{thm:quantitative2}, we have the estimates \eqref{eq:calderonestimate1}, \eqref{eq:calderonestimate2}, \eqref{eq:calderonestimate3} which imply that
        \begin{equation*}
            \left\|K_{i,j}'\right\|_{CZ} \lesssim 1,
        \end{equation*}
        with the bound independent of $0<\delta_1<\delta_2$. When combined with \eqref{eq:integralcommutator}, this is the statement of the Proposition.
    \end{proof}

    Compared with the analysis of DiPerna-Lions, in \cite{diperna1989ordinary}, the authors show that
    \begin{lemma}[Decay of the DiPerna-Lions commutator \cite{diperna1989ordinary}]\label{thm:commutatordecay}
        Let $1\le p,q\le \infty$ with $\frac{1}{p}+\frac{1}{q}=\frac{1}{r}$ with $1\le r<\infty$. Let $u(x) \in W_x^{1,p}$ be divergence-free, and $\rho(x)\in L_x^q$ then
        \begin{equation*}
            R_\delta(u,\rho)(x) \xrightarrow{\delta\to0} 0 \text{ in } L_x^r,
        \end{equation*}
    \end{lemma}
    \begin{proof}
        See \cite{diperna1989ordinary}.
    \end{proof}

    Quantitative estimates for passive scalar transport then follow from quantifying this decay. However, in the original proof \cite{diperna1989ordinary}, this decay is not uniform in the norms of $\nabla u(x) \in L_x^{p}$ or $\rho(x)\in L_x^q$. Proposition \ref{thm:expdecay} partly addresses this problem by giving a uniform integral estimate on the decay of the commutator, and the main result of this paper, Theorem \ref{thm:quantitative2}, is proved similarly. We now ask if we can improve this to quantitative estimates on the decay of the \textit{norm} of the commutator, $\left\|R_\delta(u,\rho)\right\|_{L_x^r}\to0$. To motivate that such a estimate exists, we start by observing (by a compactness argument) that there exists some uniform decay rate for the norm:
    \begin{proposition}\label{thm:compactnessargument}
        Let $1\le p,q\le \infty$ with $\frac{1}{p}+\frac{1}{q}=1$. Let $u(x) \in W_x^{1,p}(\mathbb{R}^d)$ be divergence-free, and $\rho(x)\in L_x^q(\mathbb{R}^d)$.
        
        Fix a compact subset $B \subset \mathbb{R}^d$. If $p,q>1$, then for any $\epsilon>0$, for any $\bar{\delta}>0$, there exists $\delta_\epsilon>0$ independent of $u(x)$ and $\rho(x)$ such that
        \begin{equation*}
            \inf_{\delta\in[\delta_\epsilon,\bar{\delta}]}\left\|R_\delta(u,\rho)\right\|_{L_x^1(B)} \le \epsilon \left\|\nabla u\right\|_{L_x^{p}}\left\|\rho\right\|_{L_x^q}.
        \end{equation*}
    \end{proposition}
    \begin{proof}
        Assume for the sake of contradiction that this is not the case. Then there exists an $\epsilon>0$ and a sequence $\nabla u_n(x) \in L_x^{p}$, $\rho_n(x) \in L_x^q$ for $n\in\mathbb{N}$ with $n > \bar{\delta}^{-1}$, so that
        \begin{equation*}
            \inf_{\delta\in[\frac{1}{n},\bar{\delta}]}\left\|R_\delta(u_n,\rho_n)\right\|_{L_x^1(B)} > \epsilon \left\|\nabla u\right\|_{L_x^{p}}\left\|\rho\right\|_{L_x^q}.
        \end{equation*}
        
        By scaling we may take $\left\|\nabla u_n\right\|_{L_x^{p}},\left\|\rho_n\right\|_{L_x^q}=1$ in the above. Since $p,q>1$, we may take a weakly-* converging subsequence. Without loss of generality, we take this to be the original sequence,
        \begin{gather*}
            \nabla u_n(x) \xrightharpoonup{n\to\infty} \nabla u(x) \in L_x^{p}, \\
            \rho_n(x) \xrightharpoonup{n\to\infty} \rho(x) \in L_x^q,
        \end{gather*}
        and by the Rellich-Kondrachov theorem \cite{evans2022partial}*{Section 5.7}, we have the strong convergence
        \begin{equation*}
            u_n(x) \xrightarrow{n\to\infty} u(x) \in L_\mathrm{loc}^{p},
        \end{equation*}
        and so strong convergence on compact subsets
        \begin{equation}\label{eq:commutatorcompactnessconvergence}
            R_\delta(u_n,\rho_n)(x)\xrightarrow{n\to\infty}R_\delta(u,\rho)(x)\in L_x^1(B).
        \end{equation}
        Since we assumed that
        \begin{equation*}
            \inf_{\delta\in[\frac{1}{n},\bar{\delta}]}\left\|R_\delta(u_n,\rho_n)\right\|_{L_x^1(B)} > \epsilon,
        \end{equation*}
        then, by \eqref{eq:commutatorcompactnessconvergence}, also
        \begin{equation*}
            \inf_{\delta\in[\frac{1}{n},\bar{\delta}]}\left\|R_\delta(u,\rho)\right\|_{L_x^1(B)} > \epsilon.
        \end{equation*}
        
        However, by Lemma \ref{thm:commutatordecay}, one has
        \begin{equation*}
            \inf_{\delta\in[\frac{1}{n},\bar{\delta}]}\left\|R_\delta(u,\rho)\right\|_{L_x^1(B)} \xrightarrow{n\to\infty} 0,
        \end{equation*}
        and so we reach a contradiction.
    \end{proof}

    Our goal is to quantify this compactness argument. To this end we introduce the homogeneous Besov space $\dot{B}^s_{p,q}(\mathbb{R}^d;\mathbb{R}^{d'})$ as the Schwartz distributions $\mathcal{S}'(\mathbb{R}^d;\mathbb{R}^{d'})$ for which the following homogeneous Besov norm is finite,
    \begin{equation*}
        \left\|f\right\|_{\dot{B}^s_{p,q}(\mathbb{R}^d;\mathbb{R}^{d'})} = \left(\sum_{n=-\infty}^\infty \left(2^{ns}\left\|\psi_n * f\right\|_{L^p(\mathbb{R}^d;\mathbb{R}^{d'})}\right)^{q}\right)^{\frac{1}{q}},
    \end{equation*}
    where $\psi_n \in \mathcal{S}(\mathbb{R}^d;\mathbb{R})$ is a Littlewood-Paley decomposition, given by the Fourier transform
    \begin{equation*}
        \hat{\psi}_n(\xi) = \chi(2^{-n} \xi),
    \end{equation*}
    for a fixed choice of $\chi \in C_c^\infty(\mathbb{R}^d;\mathbb{R})$ satisfying $\mathrm{supp}\;\chi = \{\xi:\frac{1}{2}\le|\xi|\le2\}$, $\chi(\xi)>0$ if $\frac{1}{2}<|\xi|<2$, and $\sum_{n=-\infty}^\infty \chi(2^{-n}\xi)=1$ for $\xi \ne 0$, as in \cite{bergh2012interpolation}*{Chapter 6}.

    Recall the continuous embedding $L_x^p \hookrightarrow \dot{B}_{p,\mathrm{max}({p,2})}^0$ for all $1<p\le\infty$, see \cite{bahouri2011fourier}*{Theorem 2.40, 2.41}.
    
    We prove the following integral estimate on the decay of the norm of the DiPerna-Lions commutator, which improves on the integral estimate in Proposition \ref{thm:expdecay} at the expense of the decay rate. We omit the time dependence for readability:
    \begin{proposition}[Integral estimate on the norm of the DiPerna-Lions commutator]\label{thm:besovdecay}
        Let $p,p',q,r\in[1,\infty]$ with $\frac{1}{p}+\frac{1}{p'}=\frac{1}{r}$.
        
        Let $\rho(x):\mathbb{R}^d\to\mathbb{R}$ with $\rho(x)\in L^{p'}(\mathbb{R}^d;\mathbb{R})\cap \dot{B}^0_{p',q}(\mathbb{R}^d;\mathbb{R})$, and $u(x):\mathbb{R}^d\to\mathbb{R}^d$ be divergence-free with $\nabla u(x)\in L^p(\mathbb{R}^d;\mathbb{R}^d\times\mathbb{R}^d)\cap \dot{B}^0_{p,q}(\mathbb{R}^d;\mathbb{R}^d\times\mathbb{R}^d)$.

        Let $\varphi(x)\in \mathcal{S}(\mathbb{R}^d;\mathbb{R})$ (a Schwartz function) be a standard mollifier. Denote the commutator by
        \begin{equation*}
            R_\delta(u,\rho)(x) = u(x)\cdot\nabla(\rho*\varphi_\delta)(x) - \nabla\cdot((u\rho)*\varphi_\delta)(x).
        \end{equation*}
        
        Then
        \begin{equation*}
            \left(\int_0^\infty \left\|R_{\delta}(u,\rho)\right\|_{L_x^r}^q\;\frac{d\delta}{\delta}\right)^{\frac{1}{q}} \le C\left(\left\|\nabla u\right\|_{\dot{B}^0_{p,q}}\left\|\rho\right\|_{L_x^{p'}} + \left\|\nabla u\right\|_{L_x^p}\left\|\rho\right\|_{\dot{B}^0_{p',q}}\right),
        \end{equation*}
        where the constant $C>0$ depends only on the dimension $d\ge1$, and the following norm of the mollifier $\varphi(x)$: $\left\|(1+|x|^2)\nabla\varphi\right\|_{L_x^1}+\left\||x|\nabla^2\varphi\right\|_{L_x^1}$. In particular, not on the parameters $1\le p,p',q,r \le \infty$.
    \end{proposition}
    \begin{proof}
        Throughout $\lesssim$ will denote less than or equal to up to a constant depending only on the dimension $d$, and the norm of the mollifier: $\left\|(1+|x|^2)\nabla\varphi\right\|_{L_x^1}+\left\||x|\nabla^2\varphi\right\|_{L_x^1}$, and in particular not on $\delta>0$.

        We first rewrite
        \begin{align*}
            R_\delta(u,\rho)(x)&=\int_{\mathbb{R}^d}\rho(y)(u(x)-u(y))\cdot\nabla\varphi_\delta(x-y)\;dy, \\
            &=\int_{\mathbb{R}^d}\rho(x-\delta h)\left(\frac{u(x)-u(x-\delta h)}{\delta}\right)\cdot\nabla\varphi(h)\;dh,
        \end{align*}
        and so by testing against $\phi(x)\in L_x^{r'}(\mathbb{R}^d;\mathbb{R})$ for $\frac{1}{r}+\frac{1}{r'}=1$ we see immediately that
        \begin{equation}\label{eqapriori1'}
            \left\|R_\delta(u,\rho)\right\|_{L_x^r} \lesssim \frac{1}{\delta}\left\|\rho\right\|_{L_x^{p'}}\left\|u\right\|_{L_x^p}.
        \end{equation}

        Next, we consider the following approximation 
        \begin{align*}
            & R_\delta(u,\rho)(x)-\sum^d_{i,j=1}\frac{\partial u_i}{\partial x_j}(x)\int_{\mathbb{R}^d}\rho(x-\delta h)h_j\frac{\partial\varphi}{\partial h_i}\;dh \\
            & \qquad = \int_{\mathbb{R}^d}\rho(x-\delta h)\sum^d_{i=1}\left(\frac{u_i(x)-u_i(x-\delta h)}{\delta} - \sum^d_{j=1}h_j\frac{\partial u_i}{\partial x_j}(x)\right)\frac{\partial\varphi}{\partial h_i}\;dh \\
            & \qquad = \int_0^1\int_{\mathbb{R}^d}\rho(x-\delta h)\sum^d_{i,j=1}\left(h_j\frac{\partial u_i}{\partial x_j}(x-t\delta h) - h_j\frac{\partial u_i}{\partial x_j}(x)\right)\frac{\partial\varphi}{\partial h_i}\;dh\;dt \\
            & \qquad = \int_0^1\int_0^1\int_{\mathbb{R}^d}\rho(x-\delta h)\sum^d_{i,j,k=1}\left(-t\delta h_j h_k\frac{\partial^2 u_i}{\partial x_j\partial x_k}(x-st\delta h)\right)\frac{\partial\varphi}{\partial h_i}\;dh\;dt\;ds,
        \end{align*}
        and so by testing against $\phi(x)\in L_x^{r'}(\mathbb{R}^d;\mathbb{R})$ we have
        \begin{equation}\label{eqapriori2'}
            \left\|R_\delta(u,\rho)-\sum^d_{i,j=1}\frac{\partial u_i}{\partial x_j}\int_{\mathbb{R}^d}\rho(x-\delta h)h_j\frac{\partial\varphi}{\partial h_i}\;dh\right\|_{L_x^r} \lesssim \delta\left\|\rho\right\|_{L_x^{p'}}\left\|\nabla^2 u\right\|_{L_x^p}.
        \end{equation}

        Now notice, since $\nabla\cdot u(x)=0$,
        \begin{align*}
            & \sum^d_{i,j=1} \frac{\partial u_i}{\partial x_j}(x)\int_{\mathbb{R}^d}\rho(x-\delta h)h_j\frac{\partial\varphi}{\partial h_i}\;dh \\
            & \qquad = \sum^d_{i,j=1}\frac{\partial u_i}{\partial x_j}(x)\int_{\mathbb{R}^d}\rho(x-\delta h)\frac{\partial (h_j \varphi)}{\partial h_i}\;dh \\
            & \qquad = \sum^d_{i,j=1}\frac{\partial u_i}{\partial x_j}(x)\int_{\mathbb{R}^d}\delta\frac{\partial \rho}{\partial x_i}(x-\delta h)h_j \varphi\;dh,
        \end{align*}
        and so by testing against $\phi(x)\in L_x^{r'}(\mathbb{R}^d;\mathbb{R})$ we have the bound
        \begin{equation}\label{eqapriori3'}
            \left\|\sum^d_{i,j=1}\frac{\partial u_i}{\partial x_j}\int_{\mathbb{R}^d}\rho(x-\delta h)h_j\frac{\partial\varphi}{\partial h_i}\;dh\right\|_{L_x^r}\lesssim\delta\left\|\nabla u\right\|_{L_x^p}\left\|\nabla \rho\right\|_{L_x^{p'}}.
        \end{equation}

        Finally, assume $\rho(x)=-\Delta g(x)$ for some $g(x) \in W_x^{1,p'}(\mathbb{R}^d;\mathbb{R})$, then
        \begin{align*}
            & \sum^d_{i,j=1}\frac{\partial u_i}{\partial x_j}\int_{\mathbb{R}^d}\rho(x-\delta h)\frac{\partial (h_j \varphi)}{\partial h_i}\;dh \\
            & \qquad = \sum^d_{i,j,k=1} \frac{\partial u_i}{\partial x_j}\int_{\mathbb{R}^d}\frac{\partial^2 g}{\partial x_k \partial x_k}(x-\delta h)\frac{\partial (h_j \varphi)}{\partial h_i}\;dh \\
            & \qquad = \sum^d_{i,j,k=1}\frac{\partial u_i}{\partial x_j}\int_{\mathbb{R}^d} -\frac{1}{\delta}\frac{\partial}{\partial h_k}\frac{\partial g}{\partial x_k}(x-\delta h)\frac{\partial (h_j \varphi)}{\partial h_i}\;dh \\
            & \qquad =\sum^d_{i,j,k=1} \frac{\partial u_i}{\partial x_j}\int_{\mathbb{R}^d}\frac{1}{\delta}\frac{\partial g}{\partial x_k}(x-\delta h)\frac{\partial^2 (h_j \varphi)}{\partial h_k \partial h_i}\;dh,
        \end{align*}
        and so
        \begin{equation}\label{eqapriori4'}
            \left\|\frac{\partial u_i}{\partial x_j}\int_{\mathbb{R}^d}\rho(x-\delta h)h_j\frac{\partial\varphi}{\partial h_i}\;dh\right\|_{L_x^r}\lesssim\frac{1}{\delta}\left\|\nabla u\right\|_{L_x^p}\left\|\nabla (-\Delta)^{-1} \rho\right\|_{L_x^{p'}}.
        \end{equation}
        
        Putting equations \eqref{eqapriori1'} to \eqref{eqapriori4'} together gives, for any $\bar{u}(x)\in W_x^{2,p}(\mathbb{R}^d;\mathbb{R}^d)$, $\bar{\rho}(x)\in W_x^{1,p'}(\mathbb{R}^d;\mathbb{R})$ such that $(-\Delta)^{-1}(\rho-\bar{\rho})(x) \in W_x^{1,p'}(\mathbb{R}^d;\mathbb{R})$,
        \begin{align}
            & \left\|R_\delta(u,\rho)\right\|_{L_x^r} \nonumber \\
            & \qquad \leq \begin{aligned}[t]
                & \left\|R_\delta(u-\bar{u},\rho)\right\|_{L_x^r} + \left\|R_\delta(\bar{u},\rho)-\sum^d_{i,j=1}\frac{\partial \bar{u}_i}{\partial x_j}\int_{\mathbb{R}^d}\rho(x-\delta h)h_j\frac{\partial\varphi}{\partial h_i}\;dh\right\|_{L_x^r} \\
                & + \left\|\sum^d_{i,j=1}\frac{\partial \bar{u}_i}{\partial x_j}\int_{\mathbb{R}^d}(\rho-\bar{\rho})(x-\delta h)h_j\frac{\partial\varphi}{\partial h_i}\;dh\right\|_{L_x^r} \\
                & + \left\|\sum^d_{i,j=1}\frac{\partial \bar{u}_i}{\partial x_j}\int_{\mathbb{R}^d}\bar{\rho}(x-\delta h)h_j\frac{\partial\varphi}{\partial h_i}\;dh\right\|_{L_x^r}
            \end{aligned} \nonumber \\
            & \qquad \lesssim \begin{aligned}[t]
                & \left\|\rho\right\|_{L_x^{p'}}\left(\frac{1}{\delta}\left\|u-\bar{u}\right\|_{L_x^p} + \delta\left\|\nabla^2\bar{u}\right\|_{L_x^p}\right) \\
                & + \left\|\nabla\bar{u}\right\|_{L_x^p}\left(\frac{1}{\delta}\left\|\nabla(-\Delta)^{-1}(\rho-\bar{\rho})\right\|_{L_x^{p'}} + \delta \left\|\nabla\bar{\rho}\right\|_{L_x^{p'}}\right).
            \end{aligned} \label{eqsplitcompact'}
        \end{align}

        We now need some standard real interpolation inequalities. For this we recall the definition of $\dot{B}^s_{p,q}$ as the Schwartz distributions $\mathcal{S}'$ such that the following norm is finite,
        \begin{equation*}
            \left\|\rho\right\|_{\dot{B}^s_{p,q}} = \left(\sum_{n=-\infty}^\infty 2^{nsq}\left\|\psi_n * \rho\right\|_{L_x^p}^q\right)^{\frac{1}{q}}.
        \end{equation*}

        By \cite{bergh2012interpolation}*{Theorem 6.4.5} the homogeneous Besov spaces $\dot{B}^1_{p,q}$, $\dot{B}^0_{p',q}$ are the real interpolation spaces between the spaces $$(\dot{B}^0_{p,1},\dot{B}^2_{p,1})_{\frac{1}{2},q} \text{ and } (\dot{B}^{-1}_{p',1},\dot{B}^1_{p',1})_{\frac{1}{2},q},$$ respectively. Thus, using the Lions-Peetre $K$-function method for real interpolation \cite{bergh2012interpolation}*{Section 3.1}, we have the following characterisation of the norms $\left\Vert \cdot \right\Vert_{\dot{B}^1_{p,q}}$ and $\left\Vert \cdot \right\Vert_{\dot{B}^0_{p',q}}$. Defining 
        \begin{gather*}
            \left| u \right|_{\dot{B}^1_{p,q}}=\bigg(\int_0^{+\infty}\inf_{\substack{\bar{u}\in \dot{B}^2_{p,1}(\mathbb{R}^d;\mathbb{R}^d)\\u-\bar{u}\in \dot{B}^{0}_{p,1}(\mathbb{R}^d;\mathbb{R}^d)}}\left(t^{-\frac{1}{2}}\left\Vert u-\bar{u}\right\Vert_{\dot{B}^0_{p,1}}+t^\frac{1}{2}\left\Vert \bar{u}\right\Vert_{\dot{B}^2_{p,1}}\right)^q\frac{dt}{t}\bigg)^{\frac{1}{q}}, \\
            \left| \rho \right|_{\dot{B}^0_{p',q}}=\bigg(\int_0^{+\infty}\inf_{\substack{\bar{\rho}\in \dot{B}^1_{p',1}(\mathbb{R}^d;\mathbb{R})\\\rho-\bar{\rho}\in \dot{B}^{-1}_{p',1}(\mathbb{R}^d;\mathbb{R})}}\left(t^{-\frac{1}{2}}\left\Vert \rho-\bar{\rho}\right\Vert_{\dot{B}^{-1}_{p',1}}+t^\frac{1}{2}\left\Vert \bar{\rho}\right\Vert_{\dot{B}^1_{p',1}}\right)^q\frac{dt}{t}\bigg)^{\frac{1}{q}},
        \end{gather*}
        there exists a universal constant $C>0$, independent of $p,p',q\in[1,\infty]$, such that 
        \begin{gather*}
            \frac{1}{C}\left\Vert u \right\Vert_{\dot{B}^1_{p,q}} \leq  \left| u \right|_{\dot{B}^1_{p,q}}\leq C \left\Vert u \right\Vert_{\dot{B}^1_{p,q}}, \\
            \frac{1}{C}\left\Vert \rho \right\Vert_{\dot{B}^{0}_{p',q}} \leq  \left| \rho \right|_{\dot{B}^0_{p',q}}\leq C \left\Vert \rho \right\Vert_{\dot{B}^0_{p',q}}.
        \end{gather*}
        Thus there exists some choice of such $\bar{u}_\delta(x), \bar{\rho}_\delta(x)$ for each $\delta>0$ so that
        \begin{gather*}
            \left(\int_0^\infty \left(\frac{1}{\delta}\left\|u-\bar{u}_\delta\right\|_{\dot{B}^0_{p,1}} + \delta\left\|\bar{u}_\delta\right\|_{\dot{B}^2_{p,1}}\right)^q \frac{d\delta}{\delta}\right)^{\frac{1}{q}} \lesssim \left\|u\right\|_{\dot{B}^1_{p,q}}, \\
            \left(\int_0^\infty \left(\frac{1}{\delta}\left\|\rho-\bar{\rho}_\delta\right\|_{\dot{B}^{-1}_{p',1}} + \delta\left\|\bar{\rho}_\delta\right\|_{\dot{B}^1_{p',1}}\right)^q \frac{d\delta}{\delta}\right)^{\frac{1}{q}} \lesssim \left\|\rho\right\|_{\dot{B}^0_{p',q}},
        \end{gather*}
        for all $q\in[1,\infty]$. Moreover, $\bar{u}_\delta(x)\in \dot{B}^2_{p,1}(\mathbb{R}^d;\mathbb{R}^d)$ may always be chosen of the form
        \begin{equation*}
            \bar{u}_\delta(x) = \sum_{n=-\infty}^{K_\delta} (\psi_n*u)(x),
        \end{equation*}
        for some $K_\delta \in \mathbb{Z}$, see \cite{bergh2012interpolation}*{Section 6.4}, and so in addition we have the bound $\left\|\nabla\bar{u}_\delta\right\|_{L_x^p}\lesssim \left\|\nabla u\right\|_{L_x^p}$ independent of $\delta>0$. In light of this, we may apply \eqref{eqsplitcompact'} for each $\delta>0$ to give the result as follows:
        \begin{align*}
            & \left(\int_0^\infty \left\|R_{\delta}(u,\rho)\right\|_{L_x^r}^q\;\frac{d\delta}{\delta}\right)^{\frac{1}{q}} \\
            & \qquad \lesssim \begin{aligned}[t]
                & \left(\int_0^\infty \left(\left\|\rho\right\|_{L_x^{p'}}\left(\frac{1}{\delta}\left\|u-\bar{u}_\delta\right\|_{L_x^p} + \delta\left\|\nabla^2\bar{u}_\delta\right\|_{L_x^p}\right)\right)^q\;\frac{d\delta}{\delta}\right)^{\frac{1}{q}} \\
                & + \left(\int_0^\infty \left(\left\|\nabla\bar{u}_\delta\right\|_{L_x^p}\left(\frac{1}{\delta}\left\|\nabla (-\Delta)^{-1}(\rho-\bar{\rho}_\delta)\right\|_{L_x^{p'}} + \delta\left\|\nabla \bar{\rho}_\delta\right\|_{L_x^{p'}}\right)\right)^q\;\frac{d\delta}{\delta}\right)^{\frac{1}{q}}
            \end{aligned} \\
            & \qquad \lesssim \begin{aligned}[t]
                & \left\|\rho\right\|_{L_x^{p'}} \left(\int_0^\infty \left(\frac{1}{\delta}\left\|u-\bar{u}_\delta\right\|_{\dot{B}^0_{p,1}} + \delta\left\|\bar{u}_\delta\right\|_{\dot{B}^2_{p,1}}\right)^q\;\frac{d\delta}{\delta}\right)^{\frac{1}{q}} \\
                & + \left\|\nabla u\right\|_{L_x^p}\left(\int_0^\infty \left(\frac{1}{\delta}\left\|\rho-\bar{\rho}_\delta\right\|_{\dot{B}^{-1}_{p',1}} + \delta\left\|\bar{\rho}_\delta\right\|_{\dot{B}^1_{p',1}}\right)^q\;\frac{d\delta}{\delta}\right)^{\frac{1}{q}}
            \end{aligned} \\
            & \qquad \lesssim \left\|\rho\right\|_{L_x^{p'}}\left\|u\right\|_{\dot{B}^1_{p,q}} + \left\|\nabla u\right\|_{L_x^p}\left\|\rho\right\|_{\dot{B}^0_{p',q}},
        \end{align*}
        where we have used \cite{bergh2012interpolation}*{Theorem 6.2.4}, that we have the continuous embedding say $\dot{B}^0_{p,1}(\mathbb{R}^d;\mathbb{R}^d) \hookrightarrow L_x^p(\mathbb{R}^d;\mathbb{R}^d)$ and we conclude the proof by $\left\|u\right\|_{\dot{B}^1_{p,q}} \lesssim \left\|\nabla u\right\|_{\dot{B}^0_{p,q}}$.
    \end{proof}

    \begin{remark}
        By the continuous embedding $L_x^p \hookrightarrow \dot{B}_{p,\mathrm{max}({p,2})}^0$ for all $1<p\le\infty$, Proposition \ref{thm:besovdecay} provides the quantitative control on $\delta_\epsilon>0$ in Proposition \ref{thm:compactnessargument}. Namely, we may take
        \begin{equation*}
            \delta_\epsilon=\bar{\delta}\exp\left(-C_{p,q}\epsilon^{-\max(p,q)}\right),
        \end{equation*}
        for some $C_{p,q}>0$ depending on the parameters $1<p,q<\infty$ and the dimension $d\ge2$.
    \end{remark}

    By the continuous embedding $L_x^p \hookrightarrow \dot{B}_{p,\mathrm{max}({p,2})}^0$ for all $1<p\le\infty$, and Minkowski's integral inequality, we have the following simple corollary:
    \begin{corollary}\label{thm:lpdecay}
        Let $\nabla u(x,t)\in L_t^pL_x^q$ be divergence-free. Let $\rho(x,t) \in L_t^{p'} L_x^{q'}$ with $\frac{1}{p}+\frac{1}{p'}=1$, and $\frac{1}{r}=\frac{1}{q}+\frac{1}{q'}\le1$.
        
        Let $\varphi(x)\in \mathcal{S}(\mathbb{R}^d;\mathbb{R})$ (a Schwartz function) be a standard mollifier. Denote the commutator by
        \begin{equation*}
            R_\delta(u,\rho)(x,t) = u(x,t)\cdot\nabla(\rho*\varphi_\delta)(x,t) - \nabla\cdot((u\rho)*\varphi_\delta)(x,t).
        \end{equation*}
        
        If $1<q,q'<\infty$, then
        \begin{equation*}
            \left(\int_0^\infty \left\|R_{\delta}(u,\rho)\right\|_{L_t^1L_x^r}^{\mathrm{max}(q,q')}\;\frac{d\delta}{\delta}\right)^{\frac{1}{{\mathrm{max}(q,q')}}} \le C_{q,q'}\left\|\nabla u\right\|_{L_t^pL_x^q}\left\|\rho\right\|_{L_t^{p'} L_x^{q'}},
        \end{equation*}
        where the constant $C_{q,q'}>0$ depends only on the parameters $1<q,q'<\infty$, the dimension $d\ge1$, and the following norm of the mollifier in $\varphi(x)$: $\left\|(1+|x|^2)\nabla\varphi\right\|_{L_x^1}+\left\||x|\nabla^2\varphi\right\|_{L_x^1}$.
    \end{corollary}
    
    Finally, we show that the integral estimates given in this section are optimal, in the sense that a non-integral estimate on the decay of the commutator does not exist (see point 1 below) and that the rate given in Proposition \ref{thm:besovdecay} is sharp when $p=p'=2$ (see point 2):
    \begin{proposition}[Sharpness of decay of the DiPerna-Lions commutator]\label{thm:counterexamples}
        Let $u(x) \in C^\infty(\mathbb{R}^2;\mathbb{R}^2)$ be the time-independent linear shear
        \begin{equation*}
            u(x) = (x_2,0) \in \mathbb{R}^2.
        \end{equation*}
        For any compactly supported mollifier $\varphi(x) \in C_c^\infty(\mathbb{R}^2;\mathbb{R})$, denote the commutator by
        \begin{equation*}
            R_\delta(u,\rho)(x) = u(x)\cdot\nabla(\rho*\varphi_\delta)(x) - \nabla\cdot((u\rho)*\varphi_\delta)(x).
        \end{equation*}
        
        Then
        \begin{enumerate}
            \item for any cutoff $\chi(x) \in C_c^\infty(\mathbb{R}^2;\mathbb{R})$ with $\chi(x)=1$ whenever $|x| \le 1$,
            \begin{equation*}
                \liminf_{\delta\to0}\sup_{\substack{\rho(x)\in C_c^\infty(\mathbb{R}^2;\mathbb{R})\\|\rho(x)|\le \chi(x)}}\left\|R_\delta(u,\rho)\right\|_{L_x^1} > 0.
            \end{equation*}
            \item for some particular choice of $\bar{\rho}(x):\mathbb{R}^2 \to \mathbb{R}$ which is in $L_x^p$ for all $p<\infty$, and supported on $|x| \le 1$, then for any $1\le q <2$
            \begin{equation*}
                \left(\int_0^1 \left\|R_\delta(u,\bar{\rho})\right\|_{L_x^1}^q \frac{d\delta}{\delta}\right)^{\frac{1}{q}} = \infty,
            \end{equation*}
        \end{enumerate}
    \end{proposition}
    \begin{proof}
        For any $\rho(x)\in L_x^1(\mathbb{R}^2;\mathbb{R})$,
        \begin{align*}
            R_\delta(u,\rho)(x) & = \int_{\mathbb{R}^2} \rho(y) (x_2-y_2)\frac{\partial \varphi_\delta}{\partial x_1}(x-y)\; dy \\
            & = \int_{\mathbb{R}^2} \rho(x-h) \frac{\partial (h_2\varphi_\delta)}{\partial h_1}(h)\; dh, \\
            & = (\rho * K_\delta)(x),
        \end{align*}
        where $K_\delta(x) = \delta^{-2} K(\delta^{-1}x)$ and $K(x) = \frac{\partial (x_2 \varphi)}{\partial x_1}(x)$.

        In particular, the Fourier transforms of $K_\delta(x)$ and $K(x)$ are given by
        \begin{equation}\label{eq:commutatorfourier}
            \begin{gathered}
                \hat{K}_\delta(\xi) = \hat{K}(\delta\xi), \\
                \hat{K}(\xi) = - \xi_1 \frac{\partial\hat{\varphi}}{\partial \xi_2}.
            \end{gathered}
        \end{equation}

        For now consider any $\rho(x)\in\mathcal{S}'(\mathbb{R}^2;\mathbb{R})$ a Schwartz distribution. Then $(\rho*K_\delta)(x)$ will coincide with $R_\delta(u,\rho)(x)$ if in addition $\rho(x)\in L_x^1(\mathbb{R}^2;\mathbb{R})$.
        
        We now prove the first statement of the theorem:
        
        Consider the harmonic $\rho(x) = e^{i\delta^{-1}\bar{\xi}\cdot x}$ for some $\bar{\xi}\in\mathbb{R}^2$ such that $|\hat{K}(\bar{\xi})|>0$, which exists since $\varphi(x)\in L_x^1(\mathbb{R}^2;\mathbb{R})$ with $\int_{\mathbb{R}^2}\varphi(x)\;dx=1$. Then $(\rho*K_\delta)(x) = e^{i\delta^{-1}\bar{\xi}\cdot x} \hat{K}(\bar{\xi})$.
        
        We now cutoff $\rho(x)$ so that it lies in $L_x^1(\mathbb{R}^d;\mathbb{R})$. Let the mollifier $\varphi(x)$ be supported on $|x| \le R$, then $\varphi_\delta(x)$ is supported on $|x| \le \delta R$. Taking a compactly supported cutoff $\chi(x)$ with $\chi(x)=1$ on $|x| \le 1$, then $R_\delta(u,\chi \rho)(x) = (\rho * K_\delta)(x)$ whenever $|x| \le 1-\delta R$. Therefore, for all $\delta \le \frac{1}{2R}$, $\left\|R_\delta(u,\chi \rho)\right\|_{L_x^1} \ge |\hat{K}(\bar{\xi})|$ as required.

        Note that $(\chi \rho)(x)$ is currently complex-valued, but by taking either the real or imaginary part, a similar non-zero bound must hold.

        We now prove the second statement of the theorem:
        
        Take some $\bar{\xi}\in\mathbb{R}^2$ with $|\bar{\xi}|=1$, and in addition some $0<c_1<c_2 <2c_1$, such that $\inf_{\delta \in [c_1,c_2]} |\hat{K}(\delta \bar{\xi})| = \epsilon >0$, noting that $\hat{K}(\xi)$ is continuous if $x_2 \varphi(x) \in L_x^1(\mathbb{R}^2;\mathbb{R})$.

        Let $\rho(x)\in \mathcal{S}'(\mathbb{R}^2;\mathbb{R})$ be of the form $\rho(x) = \sum_{n=1}^\infty a_n e^{i\delta_n^{-1}\bar{\xi}\cdot x}$ for some $a_n \in [0,1]$ and $\delta_n \in (0,\infty)$.

        Now, by \eqref{eq:commutatorfourier},
        \begin{align*}
            \left|\left(e^{i\delta_n^{-1}\bar{\xi}\cdot x}*K_\delta\right)(x)\right| & = \left|\hat{K}\left(\frac{\delta}{\delta_n}\bar{\xi}\right)\right| \\
            & \le \min\left\{\frac{\delta}{\delta_n}\left\|\frac{\partial\hat{\varphi}}{\partial\xi_2}\right\|_{L_\xi^\infty},\frac{\delta_n}{\delta}\left\||\xi|^2\frac{\partial\hat{\varphi}}{\partial\xi_2}\right\|_{L_\xi^\infty}\right\} \\
            & \le C\min\left\{\frac{\delta}{\delta_n},\frac{\delta_n}{\delta}\right\}.
        \end{align*}

        Now let $\delta_n = 2^{-n^2}$. Then if $\frac{\delta}{\delta_n} \in [c_1,c_2] \subset [c_1, 2c_1]$,
        \begin{align*}
            |(\rho*K_\delta)(x)| & \ge a_n\left|\hat{K}\left(\frac{\delta}{\delta_n}\bar{\xi}\right)\right| - C \sum_{1\le m \le n-1} |a_m|\frac{c_1 2^{-n^2}}{2^{-m^2}} - C \sum_{m = n+1}^{\infty} |a_m|\frac{2^{-m^2}}{c_2 2^{-n^2}} \\
            & \ge a_n\left|\hat{K}\left(\frac{\delta}{\delta_n}\bar{\xi}\right)\right| - c_1C2^{2-2n} - c_2^{-1}C2^{-2n} \\
            & \ge a_n\epsilon - C(c_1+c_2^{-1})4^{1-n}.
        \end{align*}

        Therefore, for any $q \in [1,2)$,
        \begin{align*}
            \int_{\delta_n c_1}^{\delta_n c_2} \left(\int_{\left[-\frac{1}{2},\frac{1}{2}\right]^2} |(\rho*K_\delta)(x)| \; dx\right)^q \frac{d\delta}{\delta} & \ge \left(\max\left\{0,a_n \epsilon - C(c_1+c_2^{-1})4^{1-n} \right\}\right)^q\log\left(\frac{c_2}{c_1}\right) \\
            & \ge \left(\frac{1}{2} a_n^q \epsilon^q - (C(c_1+c_2^{-1})4^{1-n})^q \right)\log\left(\frac{c_2}{c_1}\right),
        \end{align*}
        by the inequality $\max\{0,(a-b)\}^q + b^q \ge \frac{1}{2}a^q$ for $a,b\ge0$, $q\le2$. Then, using that $\delta_{n+1}c_1<\delta_n c_2$ are disjoint since $c_2<2c_1$, taking some $N\in\mathbb{N}$ such that $2\delta_N c_1 \le 1$,
        \begin{align*}
            \int_0^1 \left(\int_{\left[-\frac{1}{2},\frac{1}{2}\right]^2} |(\rho*K_\delta)(x)| \; dx\right)^q \frac{d\delta}{\delta} \ge \frac{1}{2}\left(\epsilon^q\sum_{n=N}^\infty a_n^q\right)\log\left(\frac{c_2}{c_1}\right) - \frac{C^q(c_1+c_2^{-1})^q}{1-4^{-q}}\log\left(\frac{c_2}{c_1}\right) .
        \end{align*}

        Now, choose $\{a_n\}_{n=1}^\infty\subset (0,\infty)$ so that $\sum_{n=N}^\infty a_n^q = \infty$ but $\sum_{n=N}^\infty a_n^2 <\infty$.
        Now, identify $\left[-\frac{1}{2},\frac{1}{2}\right]^2$ with the torus $\mathbb{T}^2$. Then $\rho(x) \in \dot{B}^0_{p,2}\left(\mathbb{T}^2;\mathbb{R}\right)$ for all $p\in[1,\infty)$ (noting that $\rho(x)$ is periodic with period 1, and each $e^{i\delta_n^{-1}\bar{\xi}\cdot x}$ lies in a different Littlewood-Paley block). By the embedding $B^0_{p,2}(\mathbb{T}^2)\hookrightarrow L_x^p(\mathbb{T}^2)$ whenever $p \in [2,\infty)$ see \cite{schmeisser1987topics}, we have that $\rho(x) \in L_x^p\left(\left[-\frac{1}{2},\frac{1}{2}\right]^2;\mathbb{R}\right)$ for all $p\in [2,\infty)$.

        We now cut $\rho(x)$ so that it lies in $L_x^1(\mathbb{R}^d;\mathbb{R})$. Let the mollifier $\varphi(x)$ be supported on $|x| \le R$, then $\varphi_\delta(x)$ is supported on $|x| \le \delta R$. Taking a compactly supported cutoff $\chi(x)$ with $\chi(x)=1$ on $|x| \le \frac{3}{4}$ and $\chi(x) = 0$ on $|x| \ge 1$, then $(\chi \rho)(x) \in L_x^p(\mathbb{R}^2)$ for all $p\in[1,\infty)$, and $R_\delta(u,\chi \rho)(x) = (\rho*K_\delta)(x)$ provided that $|x| \le \frac{3}{4}-\delta R$.

        Therefore,
        \begin{equation*}
            \int_0^1\left\|R_\delta(u,\rho\chi)\right\|_{L_x^1}^q \frac{d\delta}{\delta} = \infty.
        \end{equation*}

        Note that $(\chi \rho)(x)$ above is complex-valued, but the above integral must also be infinite for either the real or imaginary part.
    \end{proof}

\addcontentsline{toc}{section}{References}
\bibliographystyle{apa}
\bibliography{references.bib}

\end{document}